\documentclass[10pt]{article}

\usepackage[a4paper]{geometry}
\usepackage[show]{ed}

\usepackage{amsmath}
\usepackage{amsfonts}
\usepackage{latexsym}
\usepackage{amssymb}
\usepackage{graphicx}
\usepackage{subfigure} 

\usepackage{mathtools}
\usepackage{mleftright}
\usepackage{booktabs}
\usepackage{tikz} 
\usepgflibrary{shapes.misc}
\usepackage{makecell}
\usepackage{cellspace}
\usepackage[dotocloa,boxed]{algorithm2e}

\newcommand{\Rea}{\mathbb{R}}
\newcommand{\rd}{{\rm d}}
\newcommand{\domain}{D}
\newcommand{\OR}{\domain}
\newcommand{\HoDkk}{{H^1_{k}(\domain)}}

\newcommand{\CC}{\mathbb C}

\newcommand{\EE}{\mathbb E}

\newcommand{\RR}{\mathbb R}

\newcommand{\cA}{\mathcal A}

\newcommand{\cC}{\mathcal C}

\newcommand{\cK}{\mathcal K}

\newcommand{\cO}{\mathcal O}

\newcommand{\cT}{\mathcal T}

\newcommand{\bff}{\mathbf{f}}

\newcommand{\br}{\mathbf{r}}

\newcommand{\bu}{\mathbf{u}}
\newcommand{\bv}{\mathbf{v}}
\newcommand{\bw}{\mathbf{w}}
\newcommand{\bx}{\mathbf{x}}
\newcommand{\by}{\mathbf{y}}

\newcommand{\bze}{\mathbf{0}}

\newcommand{\beq}{\begin{equation}}
\newcommand{\eeq}{\end{equation}}
\newcommand{\beqs}{\begin{equation*}}
\newcommand{\eeqs}{\end{equation*}}
\newcommand{\bit}{\begin{itemize}}
\newcommand{\eit}{\end{itemize}}
\newcommand{\ben}{\begin{enumerate}}
\newcommand{\een}{\end{enumerate}}
\newcommand{\bal}{\begin{align}}
\newcommand{\eal}{\end{align}}
\newcommand{\bals}{\begin{align*}}
\newcommand{\eals}{\end{align*}}
\newcommand{\bse}{\begin{subequations}}
\newcommand{\ese}{\end{subequations}}
\newcommand{\bpr}{\begin{proposition}}
\newcommand{\epr}{\end{proposition}}
\newcommand{\bre}{\begin{remark}}
\newcommand{\ere}{\end{remark}}
\newcommand{\bpf}{\begin{proof}}
\newcommand{\epf}{\end{proof}}
\newcommand{\ble}{\begin{lemma}}
\newcommand{\ele}{\end{lemma}}
\newcommand{\bco}{\begin{corollary}}
\newcommand{\eco}{\end{corollary}}
\newcommand{\bex}{\begin{example}}
\newcommand{\eex}{\end{example}}
\newcommand{\bth}{\begin{theorem}}
\newcommand{\enth}{\end{theorem}}
\newcommand{\bcon}{\begin{condition}}
\newcommand{\econ}{\end{condition}}
\newcommand{\bas}{\begin{assumption}}
\newcommand{\eas}{\end{assumption}}
\newcommand{\bde}{\begin{definition}}
\newcommand{\ede}{\end{definition}}

\newcommand{\ton}{\text{ on }}
\newcommand{\tin}{\text{ in }}
\newcommand{\tfa}{\text{ for all }}
\newcommand{\tfor}{\text{ for }}

\newcommand{\tand}{\text{ and }}
\newcommand{\tst}{\text{ such that }}
\newcommand{\tif}{\text{ if }}

\newcommand{\tfae}{\text{ for almost every }}

\newcommand{\minispace}{\;\!}
\newcommand{\noi}{\noindent}

\newcommand{\tforall}{\text{ for all }}

\newcommand{\AmatoI}{\matI{\Amato}}
\newcommand{\Amato}{\matrixAone}

\newcommand{\Amat}{\mat{A}}

\newcommand{\mat}[1]{\mathbf{#1}}
\newcommand{\Amatt}{\matrixAtwo}
\newcommand{\matrixA}{\Amat}
\newcommand{\matrixAone}{\Amat_1}
\newcommand{\matrixAtwo}{\Amat_2}
\newcommand{\eps}{\varepsilon}
\newcommand{\grad}{\nabla}

\newcommand{\matI}[1]{\mleft(#1\mright)^{-1}}

\newcommand{\butilde}{\widetilde{\uvec}}

\newcommand{\Imat}{\mat{I}}
\newcommand{\Amin}{A_{\min}}
\newcommand{\Amax}{A_{\max}}
\newcommand{\EXP}[1]{\EE\mleft[#1\mright]}

\newcommand{\Dm}{D_{-}}
\newcommand{\Dp}{D_{+}}
\newcommand{\de}{:=}
\newcommand{\RRd}{\RR^{d}}
\newcommand{\Dmclos}{\clos{\Dm}}
\newcommand{\clos}[1]{\overline{#1}}
\newcommand{\GD}{\Gamma_{D}}

\newcommand{\Lt}[1]{L^{2}\mleft(#1\mright)}

\newcommand{\Hh}[1]{H^{1/2}\mleft(#1\mright)}

\newcommand{\Li}[1]{L^{\infty}\mleft(#1\mright)}
\newcommand{\nmin}{n_{\min}}
\newcommand{\nmax}{n_{\max}}

\newcommand{\abs}[1]{\mleft|#1\mright|}
\newcommand{\bxi}{\boldsymbol{\xi}}
\newcommand{\bxibar}{\overline{\bxi}}
\newcommand{\CCd}{\CC^{d}}

\newcommand{\defn}[1]{\emph{#1}}

\newcommand{\trace}{\gamma}

\newcommand{\dn}{\partial_{\nu}}
\newcommand{\gI}{g_{I}}
\newcommand{\Dtilde}{\widetilde{D}}
\newcommand{\compcont}{\subset\subset}
\newcommand{\GI}{\Gamma_{I}}
\newcommand{\LtD}{\Lt{D}}
\newcommand{\LtGI}{\Lt{\GI}}
\newcommand{\LiDRR}{\Li{D;\RR}}
\DeclareMathOperator{\supp}{supp}
\newcommand{\LiDRRdtd}{\Li{D;\spd}}
\newcommand{\LiDRRdtdop}{\Li{D;\op}}

\newcommand{\HoD}{\Ho{D}}
\newcommand{\Ho}[1]{H^{1}\mleft(#1\mright)}
\newcommand{\trGI}{\trace_{I}}
\newcommand{\HozD}[1]{H^{1}_{0,\mathrm{D}}\mleft(#1\mright)}

\newcommand{\vbar}{\overline{v}}

\newcommand{\DPGI}[2]{\langle #1, #2\rangle_{\GI}}
\newcommand{\DtN}{{\rm DtN}}

\newcommand{\GR}{\Gamma_{R}}

\newcommand{\N}[1]{\mleft\|#1\mright\|}
\newcommand{\Co}{C_{1}}

\newcommand{\NLtD}[1]{\N{#1}_{\LtD}}

\newcommand{\Cttilde}{\widetilde{C}_{2}}

\newcommand{\set}[1]{\mleft\{#1\mright\}}
\newcommand{\vh}{v_{h}}

\newcommand{\Vhp}{V_{h,p}}

\newcommand{\uh}{u_{h}}

\newcommand{\Th}{\cT_{h}}

\newcommand{\kz}{k_{0}}

\newcommand{\HhGI}{\Hh{\GI}}

\newcommand{\half}{\frac{1}{2}}

\newcommand{\Hmh}[1]{H^{-1/2}\mleft(#1\mright)}

\newcommand{\Amatj}{\Amat_j}

\newcommand{\nt}{n_{2}}
\newcommand{\Dmat}{\mat{D}}

\newcommand{\Hok}[1]{H^{1}_{k}\mleft(#1\mright)}

\newcommand{\NDmatk}[1]{\N{#1}_{\Dmatk}}
\newcommand{\Dmatk}{\Dmat_{k}}
\newcommand{\NDmatkI}[1]{\N{#1}_{\DmatkI}}
\newcommand{\DmatkI}{\Dmat_{k}^{-1}}
\newcommand{\vi}{v_{i}}
\newcommand{\phii}{\phi_{i}}

\newcommand{\mpm}{m_{\pm}}

\newcommand{\Smat}{\mat{S}}
\newcommand{\Mmat}{\mat{M}}
\newcommand{\Nmat}{\mat{N}}

\newcommand{\tbu}{\butilde}
\newcommand{\tu}{\utilde}
\newcommand{\utilde}{\widetilde{u}}
\newcommand{\matrixC}{\Cmat}
\newcommand{\Cmat}{\mat{C}}
\newcommand{\CCN}{\CC^{N}}

\newcommand{\re}{\mathrm{e}}

\newcommand{\LtDCC}{\Lt{D;\CC}}
\newcommand{\HmhGI}{\Hmh{\GI}}

\newcommand{\Cfn}[1]{C_{#1}}

\newcommand{\ftilde}{\widetilde{f}}

\newcommand{\HozDD}{\HozD{D}}

\newcommand{\Nfn}[2]{\N{#2}_{#1}}

\newcommand{\Cinvs}{\Cfn{\mathrm{inv},s}}
\newcommand{\NDk}[1]{\Nfn{\Dmatk}{#1}}
\newcommand{\NDkI}[1]{\Nfn{\DmatkI}{#1}}
\newcommand{\Mmatn}{\Mmat_{n}}
\newcommand{\SmatA}{\Smat_{A}}
\newcommand{\mplus}{m_{+}}
\newcommand{\mminus}{m_{-}}

\newcommand{\Nt}[1]{\Nfn{2}{#1}}
\newcommand{\splus}{s_{+}}

\newcommand{\rvecz}{\rvec_{0}}
\newcommand{\rvecm}{\rvec_{m}}

\newcommand{\cTh}{\cT_{h}}

\newcommand{\HokD}{\Hok{D}}
\newcommand{\NHok}[2]{\Nfn{\Hok{#1}}{#2}}

\newcommand{\NHokD}[1]{\NHok{D}{#1}}

\newcommand{\NLqDRR}[1]{\Nfn{\LqDRR}{#1}}

\newcommand{\LqDRR}{\Lq{D,\RR}}

\newcommand{\omegasl}{\omega_{\ell}}

\newcommand{\NLiDRR}[1]{\Nfn{\LiDRR}{#1}}

\newcommand{\HozDDs}{(\HozDD)^{*}}
\newcommand{\Ftilde}{\widetilde{F}}

\newcommand{\CCNtN}{\CC^{N\times N}}

\newcommand{\xvecz}{\xvec^{0}}

\newcommand{\NLsD}[1]{\Nfn{\LsD}{#1}}

\newcommand{\LsD}{\Lsspace{D}}
\newcommand{\Lsspace}[1]{L^{s}\mleft(#1\mright)}

\newcommand{\LqDRRdtdop}{\Lq{D,\op}}
\newcommand{\NLqDRRdtd}[1]{\Nfn{\LqDRRdtdop}{#1}}
\newcommand{\ceil}[1]{\mleft\lceil#1\mright\rceil}

\newcommand{\femvector}[1]{\mathbf{#1}}
\newcommand{\uvec}{\femvector{u}}

\newcommand{\vvec}{\femvector{v}}
\newcommand{\wvec}{\femvector{w}}
\newcommand{\fvec}{\femvector{f}}
\newcommand{\dvec}{\femvector{d}}
\newcommand{\rvec}{\femvector{r}}
\newcommand{\xvec}{\femvector{x}}

\newcommand{\aG}{a}

\newcommand{\T}{T}

\newcommand{\NLiDRRdtd}[1]{\N{#1}_{\LiDRRdtdop}}

\newcommand{\spd}{\mathsf{SPD}}

\newcommand{\coeffAj}{A_j}
\newcommand{\coeffA}{A}
\newcommand{\coeffnj}{n_j}
\newcommand{\coeffn}{n}
\newcommand{\solj}{u_j}
\newcommand{\hatx}{\widehat{\bx}}
\newcommand{\Qtwo}{Q1$^\prime$ }

\numberwithin{equation}{section}

\newtheorem{theorem}{Theorem}[section]

\newenvironment{proof}[1][Proof]{\noindent\emph{#1}\,}{\hfill$\square$}
\newtheorem{corollary}[theorem]{Corollary}
\newtheorem{definition}[theorem]{Definition}
\newtheorem{lemma}[theorem]{Lemma}
\newtheorem{proposition}[theorem]{Proposition}
\newtheorem{assumption}[theorem]{Assumption}
\newtheorem{example}[theorem]{Example}
\newtheorem{remark}[theorem]{Remark}

\newtheorem{condition}[theorem]{Condition}

\newcommand{\comp}{{\rm comp}}

\DeclareMathOperator{\esssup}{ess\,sup}

\newcommand{\Lq}[1]{L^{q}\mleft(#1\mright)}

\newcommand{\tendi}{\rightarrow \infty}

\newcommand{\pdiff}[2]{\frac{\partial #1}{\partial #2}}

\newcommand{\no}{n_{1}}

\newcommand{\diff}[2]{\frac{\rd #1}{\rd #2}}

\newcommand{\ri}{\mathrm{i}} 

\newcommand{\Unif}{\mathrm{Unif}}

\newcommand{\op}{\mathrm{op}}

\newcommand{\Cotilde}{\widetilde{C}_{1}}

\usepackage{color}
\usepackage{hyperref}
\definecolor{myblue}{rgb}{0,0,0.6}
\definecolor{darkgreen}{rgb}{0,0.5,0}
\hypersetup{colorlinks=true, linkcolor=myblue,citecolor=myblue,filecolor=myblue,urlcolor=myblue}

\definecolor{escol}{rgb}{0,0,0.8}
\definecolor{orpcol}{rgb}{0,0.6,0}

\usepackage[normalem]{ulem} 

\begin{document}

\title{Analysis of a Helmholtz preconditioning problem motivated by uncertainty quantification}
\author{I.~G.~Graham, O.~R.~Pembery, E.~A.~Spence\footnotemark[1]}

\footnotetext[1]{Department of Mathematical Sciences, University of Bath, Bath, BA2 7AY, UK, \tt I.G.Graham@bath.ac.uk, opembery@gmail.com, E.A.Spence@bath.ac.uk }

\maketitle

\begin{abstract}
This paper analyses the following question: let $\matrixA_j$, $j=1,2,$ be the Galerkin matrices corresponding to finite-element discretisations of the exterior Dirichlet problem for the heterogeneous Helmholtz equations $\nabla\cdot (A_j \nabla u_j) + k^2 n_j u_j= -f$. How small must 
$\|A_1 -A_2\|_{L^q}$ and 
$\|{n_1} - {n_2}\|_{L^q}$ be (in terms of $k$-dependence) for GMRES applied to either $(\matrixAone)^{-1}\matrixAtwo$ 
or $\matrixAtwo(\matrixAone)^{-1}$
to converge in a $k$-independent number of iterations
 for arbitrarily large $k$? (In other words, for $\matrixAone$ to be a good left- or right-preconditioner for $\matrixAtwo$?).
We prove results answering this question, give theoretical evidence for their sharpness, and give numerical experiments supporting the estimates.

Our motivation for tackling this question comes from calculating quantities of interest for the Helmholtz equation with \emph{random} coefficients $A$ and $n$. Such a calculation may require the solution of many deterministic Helmholtz problems, each with different $A$ and $n$, and the answer to the question above dictates to what extent a previously-calculated inverse of one of the Galerkin matrices can be used as a preconditioner for other Galerkin matrices. 



\end{abstract}

\section{Introduction}\label{sec:intronbpc}

\subsection{Statement of the problem}\label{sec:problem}

Let $D_- \subset \Rea^d, d=2,3, $ be a bounded Lipschitz open set such that the open complement $D_+:= \Rea^d\setminus\overline{D_-}$ is connected. 
Given $f$ with compact support and coefficients $A_j$, $n_j$, $j=1,2$,
let $\solj$, $j=1,2,$ satisfy 
the heterogeneous Helmholtz equation 
\beq\label{eq:pde}
\nabla\cdot(\coeffAj \nabla\solj ) + k^2 \coeffnj \solj =-f \quad \tin D_+,
\eeq
the Dirichlet boundary condition $u_j=0$ on $\Gamma:= \partial D_-$, and the Sommerfeld radiation condition
\beq\label{eq:src}
\pdiff{u_j}{r}(\bx) - \ri k u_j(\bx) = o \left( \frac{1}{r^{(d-1)/2}}\right)
\eeq
as $r:= \|\bx\|_2\tendi$ (uniformly in $\hatx:= \bx/r$),
where $k>0$ is the wavenumber; in this paper, we are interested in the case when the wavenumber is arbitrarily large.
In Equation \eqref{eq:pde} the obstacle $D_-$ and the coefficients $A_j$ and $n_j$ satisfy the natural conditions for the scattering problem to make sense; see Assumption \ref{ass:1} and Definition \ref{def:TEDP}
(to present the main results as close to the beginning of the paper as possible, we postpone to \S\ref{sec:BVP}
the precise definitions of the Helmholtz problem, the finite-element method, and GMRES).

Let $\matrixA_j$, $j=1,2,$ be the Galerkin matrices corresponding to the $h$-version finite-element discretisations (with decreasing mesh size $h$ and any fixed polynomial degree $p$)
of \eqref{eq:pde} truncated to ${D}:= D_+ \cap \widetilde{D}$, where $\widetilde{D}$ is some open set containing $\Dm$
(the simplest case is when $\widetilde{D}= B_R:= \{ \bx : \|\bx\|_2<R\}$ for some sufficiently large $R$).
We consider the cases where  \emph{either} the radiation condition is realised exactly on $\partial \widetilde{D}$
via the exact Dirichlet-to-Neumann (DtN) map \emph{or} the DtN map is approximated on $\partial \widetilde{D}$ by an impedance boundary condition (a.k.a.~a first-order absorbing boundary condition).
See \eqref{eq:matrixAjdef} for the precise definitions of $\matrixA_j$, $j=1,2$, and Definition \ref{def:TEDP} and Remark \ref{rem:T} for the truncation procedure.
 
This paper considers the following question:

\bit
\item[Q1.] How small must $\|\coeffA_1 - \coeffA_2\|_{L^q}$ and 
$\|n_1 - n_2\|_{L^q}$ be (in terms of $k$-dependence) for the generalised minimum residual method (GMRES) applied to either $(\matrixAone)^{-1}\matrixAtwo$
or  $(\matrixAtwo)^{-1}\matrixAone$
 to converge in a $k$-independent number of iterations
 for arbitrarily large $k$? 
\eit
Before stating results answering Q1 in \S\ref{sec:mainresults} we describe our motivation for studying Q1.

\subsection{Motivation from uncertainty quantification of the Helmholtz equation}\label{sec:UQ1}
Consider the Helmholtz equation 
\beq\label{eq:pde_stochastic}
\nabla\cdot\big(\coeffA(\bx;\omega) \nabla u(\bx,\omega) \big) + k^2 \coeffn(\bx;\omega) u(\bx;\omega) =-f(\bx), \quad \bx\in D_+,
\eeq
where $\coeffA(\bx;\omega)$ and $\coeffn(\bx;\omega)$ are \emph{random fields}, and $\omega$ is an element of an underlying probability space.
Suppose $Q(u)$ is a functional of the solution of \eqref{eq:pde_stochastic} (usually called a quantity of interest). In the forward problem of uncertainty quantification (UQ), one task is to compute $\EXP{Q(u)}$,
and the arguably 
simplest way to do this uses sampling, i.e.~using the approximation
\beq\label{eq:samplingexp}
\EXP{Q(u)} \approx \frac1N \sum_{\ell=1}^N Q(u(\omegasl)),
\eeq
where the $\omegasl$ are elements of the sample space $\Omega$. To calculate the right-hand side of \eqref{eq:samplingexp},  one must solve many deterministic Helmholtz problems, corresponding to different samples $\omegasl$, i.e. corresponding to different realisations of the coefficients $A(\cdot,\omega)$ and $n(\cdot,\omega)$.
Solving all these deterministic problems is a very computationally-intensive task because linear systems arising from discretisations of the Helmholtz equation are notoriously difficult to solve;
this difficulty is due to the following three reasons:
\begin{enumerate}
\item[(a)] 
The solutions of the homogeneous Helmholtz equation $\Delta u +k^2 u=0$ oscillate on a scale of $1/k$, and so to approximate them accurately with piecewise-polynomial functions one needs the total number of degrees of freedom, $N$, to be proportional to $k^d$ as $k$ increases.
\item[(b)] The \emph{pollution effect} means that, for fixed-order finite-element methods, having $N\sim k^d$ is still not enough to keep the relative finite-element error bounded independently of $k$ as $k$ increases (see the references in \S\ref{sec:cond2hold} for how $N$ must depend on $k$).
\item[(c)]
When iterative methods such as GMRES are applied to the linear system $\matrixA \uvec=\fvec$, the number of iterations grows with $k$.
Part of the explanation for this is that the standard variational formulation of the Helmholtz equation is not coercive (i.e.~it is sign-indefinite) when $k$ is sufficiently large,
and this indefiniteness is inherited by the Galerkin linear system. The design of good preconditioners for discretisations of the Helmholtz equation is therefore a very active area of research; see, e.g., the literature reviews in \cite[\S1.3]{GrSpZo:18}, \cite[\S4]{GaZh:19}.
\end{enumerate}

The question therefore arises of how to compute $Q(u(\omegasl))$, $\ell=1,\ldots,N$, in an efficient way. 
To simplify the discussion, initially suppose that one calculates the LU factorisation of the Galerkin matrix for one realisation of the coefficients $A$ and $n$. The answer to Q1 above makes $k$-explicit the extent to which this factorisation can be used as an effective preconditioner for Galerkin matrices arising from different realisations of $A$ and $n$.

The benefit of ``reusing" an LU factorisation in this way can be seen by recalling that direct solvers involving the LU decomposition of the linear system have a computational cost of the order $\cO\mleft(N^{3/2}\mright)$ in 2-d \cite[Section 1]{DuErRe:76}  and $\cO\mleft(N^2\mright)$ in 3-d \cite[Equation 3]{DuErRe:76}, with this analysis assuming sufficient grid regularity. 
However, if one already has computed the LU decomposition, then the cost of applying backsolves using it 
is much cheaper; $\cO\mleft(N\log N\mright)$ in 2-d \cite[Section 1]{DuErRe:76} and $\cO\mleft(N^{4/3}\mright)$ in 3-d \cite[Equation 4]{DuErRe:76}, hence the LU decomposition of $\Amat$ could be efficiently used as a preconditioner for matrices ``near" $\Amat$. 

More generally, for \emph{any} preconditioner for the Helmholtz equation, one could ask when the preconditioner corresponding to one realisation of \eqref{eq:pde_stochastic} can be re-used for other realisations. The answer to Q1 does not answer this more-complicated scenario, but forms the first step in analysing it. 

In \S\ref{sec:UQ2} we discuss the implications of our answer to Q1 on such a preconditioning strategy. However, 
we highlight that to perform a complete $k$-explicit analysis of this type of preconditioning strategy applied to computing $\EXP{Q(u)}$, 
in addition to the answer to Q1, one also needs the answers to the following two questions.
\bit
\item[Q2] How must the number of samples, $N$, depend on $k$ for the error
\beq\label{eq:QMCerror}
\left|\EXP{Q(u)} -\frac1N \sum_{\ell=1}^N Q(u(\omegasl))\right|
\eeq
to be controllably small, independent of $k$, for arbitrarily large $k$?
\item[Q3] How do $\|A(\cdot,\omega_{\ell_1}) - A(\cdot,\omega_{\ell_2})\|$ and 
$\|n(\cdot,\omega_{\ell_1}) - n(\cdot,\omega_{\ell_2})\|$ depend on $\omega_{\ell_1}, \omega_{\ell_2}$?
In standard UQ applications, the random fields are parametrised in terms of vectors of i.i.d.~random variables $\by(\omega)$. In this context $A(\cdot,\omega)$, $n(\cdot,\omega)$ may be written as $A(\cdot, \by(\omega))$, $n(\cdot, \by(\omega))$; this question then asks 
how ``closeness" in parameter space (measured by $\|\by(\omega_{\ell_1}) -\by(\omega_{\ell_2})\|$) translates to ``closeness" in coefficient space.
\eit
We do not address these two questions in this paper, but remark that 
(i) the recent paper \cite{GaKuSl:20} gives results related to Q2 for a non-standard variational formulation of the heterogeneous Helmholtz equation given in \cite{GaMo:19} (based on the formulation for the homogeneous Helmholtz equation introduced in \cite{MoSp:14}), and
(ii) the answer to Q3 follows from the specific structure of the randomness assumed in the coefficients $A$ and $n$.

\subsection{Novelty of the main results}\label{sec:novelty}

To our knowledge, this paper is the first time that preconditioning a discretisation of the heterogeneous Helmholtz equation with a discretisation corresponding to a Helmholtz problem with ``nearby" coefficients has been analysed.

This paper is part of a growing body of work on the analysis and numerical analysis of the heterogeneous Helmholtz equation, with the analysis explicit in both the wavenumber and the coefficients \cite{ChNi:19, GaSpWu:20, GoGrSp:20, GrPeSp:19, GrSa:20, LaSpWu:19a, MoSp:19, Pe:20, SaTo:18}. Indeed, our main results 
(Theorems  \ref{thm:Linfty1}, \ref{thm:Linfty2}, \ref{thm:1new}, and  \ref{thm:1anew})
are proved under Conditions \ref{cond:1} and \ref{cond:2}, and these other works give sufficient conditions for 
Conditions \ref{cond:1} and \ref{cond:2} to hold; 
see \S\ref{sec:cond1hold} and \S\ref{sec:cond2hold}.

Furthermore, we believe that the idea of preconditioning with ``nearby" coefficients will be relevant to other PDEs with random coefficients (this technique 
can then then be seen as a generalisation of the mean-based preconditioning discussed in \S\ref{sec:UQ2}). We expect that the ideas in our analysis -- in particular Lemmas \ref{lem:keylemma1} and \ref{lem:keylemma2} which combine bounds on the continuous problem and on the finite-element error to prove bounds on preconditioned matrices -- would then be useful in analysing such strategies.

\paragraph{Outline of the paper.}
\S\ref{sec:mainresults} states the main results and gives supporting numerical experiments.
\S\ref{sec:BVP} gives the precise definitions of the Helmholtz boundary value problem, its finite-element solution, and weighted GMRES.
\S\ref{sec:conditions} give the precise definitions of the conditions under which the main results hold.
\S\ref{sec:proofs} and \S\ref{sec:proofsPDE} prove the main results.

\section{Statement of the main results}\label{sec:mainresults}

We state in \S\ref{sec:Linfty} the results obtained when the differences $A_1-A_2$ and $n_1-n_2$ are measured in the $L^\infty$ norm (Theorems \ref{thm:Linfty1} and \ref{thm:Linfty2}). We discuss in \S\ref{sec:sharpness} the sharpness of these results.
We then state in \S\ref{sec:Lq} generalisations of Theorems \ref{thm:Linfty1} and \ref{thm:Linfty2} where the differences are measured in the $L^q$ norm for $2< q\leq\infty$ (Theorems  \ref{thm:1new} and  \ref{thm:1anew}). We discuss in \S\ref{sec:UQ2} the implications of these results for calculating quantities of interest of the Helmholtz equation with random coefficients (following on from the discussion in \S\ref{sec:UQ1}).

We consider both standard GMRES (which uses the Euclidean norm $\|\cdot\|_2$ on vectors) and weighted GMRES,
a short description of which is given in \S\ref{sec:BVP}.
In weighted GMRES, we use the weighted vector norms $\NDmatk{\cdot}$ and $\NDmatkI{\cdot}$ defined by 
\beq\label{eq:Dk}
\NDmatk{\vvec}^2\de \big( \Dmatk \vvec, \vvec\big)_2 
\quad \tand
\quad \NDmatkI{\vvec}^2\de \big( \Dmatk^{-1} \vvec, \vvec\big)_2 
\eeq
where
\beqs
\Dmat_k\de \Smat_I + k^2 \Mmat_1,
\eeqs
with $\Smat_I$ is the standard stiffness matrix for the finite-element discretisation of the Laplacian, and $\Mmat_1$ is the mass matrix -- see \eqref{eq:matrixSjdef} for a precise definition. The key point is that, for a finite-element function $\vh =\sum_i \vi \phii$,
\beq\label{eq:Dk2}
\NDmatk{\vvec}^2 =\N{v_h}^2_{\HokD}:=  \N{\grad v_h}^2_{L^2(D)} + k^2 \N{v_h}^2_{L^2(D)}
\eeq
where $\vvec$ is the vector with $i$th entry $v_i$; i.e., $\NDmatk{\cdot}$ is the norm on the finite-element space induced by the weighted $H^1$ norm $\NHokD{\cdot}$ in which the PDE analysis of the Helmholtz equation naturally takes place.
(Results about convergence of domain-decomposition methods in the norms \eqref{eq:Dk} were recently obtained for the Helmholtz equation in  \cite{GrSpVa:17,GrSpZo:18, GoGrSp:20} and for the time-harmonic Maxwell equations in \cite{BoDoGrSpTo:19}.)
We further define $ \NLiDRRdtd{A}:=\esssup_{\bx\in D}\|A(\bx)\|_2$, where here $\|\cdot\|_2$ denotes 
the \emph{spectral/operator} norm on matrices, induced by the Euclidean norm $\|\cdot\|_2$ on vectors.

To make the statements of the main results more concise, we make the following definition for weighted/standard GMRES.

\begin{definition}\label{def:GMRESconverge}
We say that GMRES applied to $\matrixC \uvec= \fvec$, with zero starting guess and in a norm $\|\cdot\|$ on $\mathbb{C}^n$,
converges in a $k$-independent number of iterations if, given $\eps>0$ and $k_0>0$, there exists $C_1(\eps, k_0)>0$ and $C_2(k_0)>0$, both independent of $k,h,$ and $p$, such that 
\beqs
\text{ if } \, m\geq C_1 \,\, \text{ then } \,\, \frac{\N{\uvec^m - \uvec}}{\N{\uvec}}\leq C_2 \eps , \,\, \tfa k\geq k_0.
\eeqs
\end{definition}
In this paper we only consider GMRES in one of the three norms $\|\cdot\|, \|\cdot\|_{\Dmat_k},$ and $\|\cdot\|_{\Dmat^{-1}_k}$ (with the last two defined in \eqref{eq:Dk}).

\subsection{Results on Q1 involving $\|A_1-A_2\|_{L^\infty}$ and $\|n_1-n_2\|_{L^\infty}$}\label{sec:Linfty}

These results are proved under Conditions \ref{cond:1} and \ref{cond:2}. These conditions can be informally stated, respectively, as 
\bit
\item the obstacle $\Dm$ and the coefficients $A_j$ and $n_j$ are such that $u_j$ exists, is unique, and the problem is \emph{nontrapping} (in the sense described in \S\ref{sec:cond1hold}), and
\item the meshsize $h$ and polynomial degree $p$ in the finite-element method are chosen to depend on $k$ to ensure that the 
finite-element approximation to the solution of the problem with coefficients $A_j$ and $n_j$ exists, is unique, and has bounded error in the $H^1_k$-norm (defined in \eqref{eq:1knorm}) as $k\tendi$. 
\eit  

\begin{theorem}[Answer to Q1: $k$-independent weighted GMRES iterations]\label{thm:Linfty1}

\noindent Assume that 
\bit
\item $\Dm$, $A_j$, and $n_j$, $j=1,2,$ satisfy Condition \ref{cond:1}, and 
\item $\Dm$, $A_j$, $n_j$, $j=1,2$, $h$ and $p$, satisfy Condition \ref{cond:2}.
\eit
Given $k_0>0$,
there exist $C_1, C_2>0$ independent of $h$ and $k$ (but dependent on $\Dm, A_1, n_1, p$, and $\kz$) 
such that if
\beq\label{eq:cond}
C_1 \,k \,\big\|A_1-A_2\big\|_{L^\infty(D; {\rm op})} 
+C_2 \, k\, \big\|n_1-n_2\big\|_{L^\infty(D;\Rea)}
\leq \frac{1}{2}
\eeq
for all $k\geq k_0$, then \emph{both} weighted GMRES working in the vector norm $\|\cdot\|_{\Dmat_k}$ (and the associated inner product) applied to 
\beq\label{eq:pcsystem1}
(\Amat_1)^{-1}\Amat_2\bu = (\Amat_1)^{-1}\bff
\eeq
\emph{and} weighted GMRES working in the vector norm $\|\cdot\|_{(\Dmat_k)^{-1}}$ (and the associated inner product) applied to 
\beq\label{eq:pcsystem2}
\Amat_2(\Amat_1)^{-1}\bv = \bff
\eeq
with initial guess the zero vector
converge in a $k$-independent number of iterations (in the sense of Definition \ref{def:GMRESconverge}) for all $k\geq k_0$.
\end{theorem}

\begin{theorem}[Answer to Q1: $k$-independent standard GMRES iterations]\label{thm:Linfty2}

\noindent Assume that 
\bit
\item $\Dm$, $A_j$, and $n_j$, $j=1,2,$ satisfy Condition \ref{cond:1}, and 
\item $\Dm$, $A_j$, $n_j$, $j=1,2$, $h$ and $p$, satisfy Condition \ref{cond:2}.
\eit
Given $k_0>0$,
let $C_1$ and $C_2$ be as in Theorem \ref{thm:Linfty1}, and let $s_{+}$ and $m_{\pm}$ be as in Lemma \ref{lem:normequiv} (note that all of these constants are independent of $k$ and $h$). 
Then if 
\beq\label{eq:conda}
 C_1 \,\left(\frac{s_+}{m_-}\right) \,\frac{1}{h} \,
\big\|A_1-A_2\big\|_{L^\infty(D;{\rm op})} + C_2 \, \left(\frac{m_+}{m_-} \right)k \, \big\|{n_1}-{n_2}\big\|_{L^\infty(D;\Rea)}
\leq \frac{1}{2}
\eeq
for all $k\geq k_0$, then standard GMRES (working in the Euclidean norm and inner product) applied to either of the equations \eqref{eq:pcsystem1} or \eqref{eq:pcsystem2}
with initial guess the zero vector
converge in a $k$-independent number of iterations (in the sense of Definition \ref{def:GMRESconverge}) for all $k\geq k_0$.
\end{theorem}

\noi We make the following remarks.
\bit[(a)]
\item The constants $C_1$ and $C_2$ are given explicitly  in \eqref{eq:C1} and \eqref{eq:C2} in terms of 
\bit
\item $C^{(1)}_{\rm bound}$ given in Condition \ref{cond:1},
\item $C^{(1)}_{\rm FEM1}$ and $C^{(1)}_{\rm FEM2}$ given in Condition \ref{cond:2},
\item ${n_1}_{\max}$, ${n_1}_{\min}$, and ${A_1}_{\min}$ given in \eqref{eq:nbounds} and \eqref{eq:Abounds} (with $n$ replaced by ${n_1}$ and $A$ replaced by ${A_1}$),
\eit
\item[(b)] The constant $1/2$ on the right-hand sides of \eqref{eq:cond} and \eqref{eq:conda} can be replaced by any number $<1$ and the 
overall 
result that GMRES converges in a $k$-independent number of iterations still holds, although the actual number of GMRES iterations depends on this number (but not on $k$).
\item[(c)] If $A_1=A_2$, then only Part 1 of Condition \ref{cond:2} is needed. If $n_1=n_2$, then only Part 2 of Condition \ref{cond:2} is needed.
\item[(d)] When $h\sim  k^{-1}$, the bounds \eqref{eq:cond} and \eqref{eq:conda} are identical in their $k$-dependence; however, when $h\ll k^{-1}$ (as one needs to take to overcome the pollution effect, as discussed in \S\ref{sec:cond2hold}) 
and $A_1\neq A_2$, 
the bound \eqref{eq:conda} for standard GMRES is more pessimistic than the bound \eqref{eq:cond} for weighted GMRES.
\item[(e)] 
Because of the discrepancy when $A_1\neq A_2$ between our bound for standard GMRES and our bound for weighted GMRES, we include the following lemma relating standard and weighted GMRES.
This result is based on \cite[Corollary 5.8, Part (ii)]{GoGrSp:20}, and shows that if the condition \eqref{eq:cond} for $k$-independent convergence of weighted GMRES holds, then standard GMRES takes at most $\sim \log k$ additional iterations to ensure the same error 
guaranteed for weighted GMRES by Theorem \ref{thm:Linfty1}. 
\eit

\begin{lemma}[Bound on standard GMRES from condition on weighted GMRES]
\label{lem:standardvsweighted}
Given $k_0>0$, assume that \eqref{eq:cond} holds (so that weighted GMRES with initial guess the zero vector converges in a $k$-independent number of iterations). 
If $h = C k^{-1-\delta}$ with $\delta \geq 0$, then there exists $\mathcal{C}>0$, depending only on $s_+, m_\pm, C, \delta,$ and $k_0$, such that if standard GMRES is applied to the linear system \eqref{eq:pcsystem1} with initial guess the zero vector, and\beq\label{eq:compare1a}
m \geq 
\left(\log\left(\frac{3}{2\sqrt{2}}\right)\right)^{-1}
\left( \log\left(\frac{1}{\eps}\right) +  \gamma \log k + \mathcal{C}\right),
\eeq
then 
\beq\label{eq:compare1b}
\frac{\N{\uvec^m - \uvec}_2}{\N{\uvec}_2}\leq \frac{9}{4} \eps , \,\, \tfa k\geq k_0.
\eeq
\end{lemma}

\

This result is not of the same form as Theorems \ref{thm:Linfty1} and \ref{thm:Linfty2}, which give sufficient conditions for weighted/standard GMRES to converge in a $k$-independent number of iterations (in the sense of Definition \ref{def:GMRESconverge}). To better understand Lemma \ref{lem:standardvsweighted}, we observe that the proof of Theorem \ref{thm:Linfty1} (see \S\ref{sec:proofLinfty1}) shows that, given $k_0>0$, if \eqref{eq:cond} holds, and
\beq\label{eq:compare2}
\tif m\geq 
\left(\log\left(\frac{3}{2\sqrt{2}}\right)\right)^{-1} \log\left(\frac{1}{\eps}\right),
\quad\text{ then } \quad \frac{\N{\uvec^m - \uvec}_{\Dmat_k}}{\N{\uvec}_{\Dmat_k}}\leq \frac{9}{4} \eps  \qquad \tfa k\geq k_0.
\eeq
Comparing \eqref{eq:compare1a}-\eqref{eq:compare1b} and \eqref{eq:compare2}, 
we see that Lemma \ref{lem:standardvsweighted} guarantees the same error for standard GMRES as Theorem \ref{thm:Linfty1} guarantees for weighted GMRES, with the conditions on $m$ showing at most $\sim\log k$ additional iterations are required for standard GMRES.
This result indicates (but does not prove) that the factor $h^{-1}$ on the left-hand side of \eqref{eq:conda} is suboptimal, and we expect the result of Theorem \ref{thm:Linfty2} holds if this factor is replaced by $k$ (as in \eqref{eq:cond}).

\bre[Sufficient conditions for the relative residual to be controllably small]

\

\noi Theorems \ref{thm:Linfty1} and \ref{thm:Linfty2} give sufficient conditions for GMRES to converge in a $k$-independent number of iterations in the sense of Definition \ref{def:GMRESconverge}. These conditions are also sufficient for the GMRES relative residual to be controllably small, independent of $k$, in the sense that, given $\eps>0$ and $k_0>0$, there exists $C_1(\eps, k_0)>0$ and $C_2(k_0)>0$, both independent of $k,h,$ and $p$, such that 
\beq\label{eq:relresidual}
\text{ if } \, m\geq C_1 \,\, \text{ then } \,\, \frac{\N{\rvec^m}}{\N{\rvec^0}}\leq C_2 \eps , \,\, \tfa k\geq k_0.
\eeq
Indeed, the proofs of Theorems \ref{thm:Linfty1} and \ref{thm:Linfty2} first obtain conditions for \eqref{eq:relresidual} to hold via the 
Elman estimate (see Theorem \ref{thm:GMRES1_intro} below), and then bound $\|\uvec^m-\uvec\|/\|\uvec\|$ by a $k$-independent multiple of $\|\rvec^m\|/\|\rvec^0\|$.
\ere

\subsection{How sharp are Theorems \ref{thm:Linfty1} and \ref{thm:Linfty2} in their $k$-dependence?}\label{sec:sharpness}

\subsubsection{Sharpness of an analogous result at the PDE level}

We now state results showing
 that the condition
\beq\label{eq:sufficientlysmall}
k\,
\NLiDRRdtd{A_1-A_2} \quad\text{ and } \quad k\,\NLiDRR{n_1-n_2}
\quad\text{ both sufficiently small}
\eeq
is not only an answer to Q1 but is also the natural answer to an analogue of Q1 at the level of PDEs, namely 
\ben
\item[\Qtwo]
Let $u_j$, $j=1,2,$ be the Helmholtz solution with coefficients $A_j$ and $n_j$;
how small must $\NLiDRRdtd{A_1 - A_2}$ and 
$\NLiDRR{n_1 - n_2}$ be (in terms of $k$) for the relative error in approximating 
$u_2$ by $u_1$  to be bounded independently of $k$ for arbitrarily-large $k$? 
\een
(The claim that \Qtwo is an analogous question to Q1 at the level of PDEs is made precise in Remark \ref{rem:analogue}.)

Furthermore Lemma \ref{lem:sharp} shows that the condition ``$k\NLiDRR{n_1 - n_2}$ sufficiently small" is the \emph{provably-sharp} answer to \Qtwo when $A_1= A_2= I$.

Before stating these PDE results, we recall from \eqref{eq:Dk2} the definition of the weighted $H^1$ norm
\beq\label{eq:1knorm}
\N{v}^2_{\HokD} \de \N{\grad v}^2_{L^2(D)} + k^2 \N{v}^2_{L^2(D)} \quad \tfor v \in H^1_{0,D}(D),
\eeq
where the space $H^1_{0,D}(D)$, defined by \eqref{eq:spaceEDP}, is the natural space containing the solution of the exterior Dirichlet problem. 
Finally, let $(H^1_{0,D}(D))^*$ denote the dual space of $H^1_{0,D}(D)$.

\begin{theorem}[Answer to \Qtwo]\label{thm:2}
Given $F \in (H^1_{0,D}(D))^*$ and coefficients $A_j, n_j$, $j=1,2,$ let $u_j$, $j=1,2,$ be the solutions of the variational problem \eqref{eq:vp} with coefficients $A_j, n_j$.
Assume that $D$, $A_1$, and $n_1$ satisfy Condition \ref{cond:1} and that 
$D$, $A_2$, and $n_2$ are such that $u_2$ exists.

Then, given $k_0>0$, there exists $C_3>0$, independent of $k$ and $F$ and given explicitly in terms of $\Dm$, $A_1$, and $n_1$ in \eqref{eq:C3}, such that
\beq\label{eq:PDEbound}
\frac{\big\|u_1-u_2\big\|_{\HokD}
}{
\N{u_2}_{\HokD}
}\leq C_3 \,k\, \max\set{\NLiDRRdtd{A_1-A_2}\,,\, \NLiDRR{n_1-n_2}}
\eeq
for all $k\geq k_0$. 
\end{theorem}

In stating the next result (and elsewhere in the paper), 
for $a,b>0$ we write $a\lesssim b$ when $a\leq C b$ for some $C>0$, independent of $k$ and $h$. We also write $a\sim b$ if $a\lesssim b$ and $b\lesssim a$.

\ble[Sharpness of the bound \eqref{eq:PDEbound} when $A_1 = A_2= I$]\label{lem:sharp}
Assume that there exist
 $0<R_1<R_2$ such that $\overline{\Dm} \subset B_{R_1}\subset B_{R_2} \subset \widetilde{D}$ (i.e.~there exists an annulus between $\Dm$ and $\partial\widetilde{D}$).
 
Then there exist particular choices of  $F, \,n_1$, and $n_2$ (with $n_1\neq n_2$, and both continuous functions on $\Dp$) such that
the corresponding solutions $u_1$ and $u_2$ of the variational problem \eqref{eq:vp} 
with $A_1 = A_2= I$ exist, are unique, and, given $k_0>0$,
\beq\label{eq:sharp1}
\frac{\N{u_1-u_2}_{\HokD}
}{
\N{u_2}_{\HokD}
}
\sim 
\frac{\N{u_1-u_2}_{L^2(D)}
}{
\N{u_2}_{L^2(D)}
}\sim k \NLiDRR{n_1-n_2},
\eeq
for all $k\geq k_0$.
\ele

\bre[Physical interpretation of the condition \eqref{eq:sufficientlysmall}]\label{rem:physical1k}
Recall that the wavelength of the wave $u$  (at least when $A=I$ and $n=1$) is $2\pi/k$.
Since this is the natural length scale associated with $u$, one expects from physical considerations that perturbations of magnitude $\leq c/k$, for $c>0$ sufficiently small, are `unseen' by the PDE or numerical method. This is indeed what we see in Theorems \ref{thm:Linfty1}, \ref{thm:Linfty2}, and \ref{thm:2}: perturbations of this form (i.e.~$A_j$ and $n_j$ satisfying \eqref{eq:sufficientlysmall})
give bounded relative difference (for \Qtwo\hspace{-1ex}) and bounded GMRES iterations for the nearby-preconditioned linear system (for Q1). 
\ere

\subsubsection{Numerical experiments to test the sharpness of Theorems \ref{thm:Linfty1} and \ref{thm:Linfty2}}\label{sec:num}

We now present experiments testing the sharpness of the condition
\beq\label{eq:suffcondition}
k \N{A_1-A_2}_{L^\infty(D; {\rm op})} + k \N{n_1-n_2}_{L^\infty(D; {\rm op})} \text{ sufficiently small }
\eeq
for the $k$-independent convergence of standard GMRES. While Theorem \ref{thm:Linfty2} gives 

\noi ``$h \N{A_1-A_2}_{L^\infty(D; {\rm op})} + k \N{n_1-n_2}_{L^\infty(D; {\rm op})}$ sufficiently small'' as a sufficient condition for the $k$-independent convergence of standard GMRES, we recall from the discussion below Lemma \ref{lem:standardvsweighted} that we expect \eqref{eq:suffcondition} to be sufficient.
While we do not present any experiments with weighted GMRES in this paper, we highlight that in \cite{GrSpVa:17} the number of iterations for standard GMRES and weighted GMRES with the weight $\Dmat_k$ were found to be almost identical when applied to preconditioned linear systems arising from the Helmholtz FEM.

\paragraph{Description of the set-up.}

We consider the solution of the variational problem \eqref{eq:vp} with $\Dm = \emptyset$, $\Dp = [0,1]^2$, and $T=\ri k$. These choices correspond to the interior impedance problem on the unit square in 2-d. We set $A_1=I$, $n_1=1$, $f=0$, and define $\gI$ so that the exact solution $u_1$ is a plane wave incident from the bottom left passing through the homogeneous medium given by coefficients $A_1$ and $n_1$. 

We solve this variational problem via the finite-element method defined in \eqref{eq:Galerkin}; we choose the family of finite-dimensional subspaces $(V_{h,p})_{h>0}$ to be first-order (i.e.~$p=1$) continuous finite elements on regular grids with $h = k^{-3/2}$; this choice of $h$ means that the finite-element approximation $u_{h}$ is uniformly accurate as $k\tendi$ (see \S\ref{sec:cond2hold} and the references therein). All finite-element calculations were carried out using the Firedrake software library \cite{RaHaMiLaLuMcBeMaKe:16,LuVaRaBeRaHaKe:15}, which itself uses the PETSc \cite{BaGrMcSm:97,DaPaKlCo:11}, Chaco \cite{HeLe:95}, and MUMPS \cite{AmDuLEKo:01,AmGuLEPr:06} software packages. The code used in this paper can be found at \url{https://github.com/orpembery/running-nbpc}. The exact versions of the code used in this paper, and the corresponding computational data can be found in the releases \cite{Pe:20a,Pe:20c,Pe:20b}.

    \begin{figure}[h!]
    \subfigure[ $\beta = 0,0.1,0.2,0.3.$]
{
          \scalebox{0.55}{
\input{nbpc-linfinity-plot-n-0.pgf}
}
}
\subfigure[$\beta = 0.4,0.5,0.6,0.7.$]
{
          \scalebox{0.55}{
\input{nbpc-linfinity-plot-n-1.pgf}
}
}
 \subfigure[$\beta = 0.8,0.9,1.$]
{
\centering{
          \scalebox{0.55}{
\input{nbpc-linfinity-plot-n-2.pgf}
}
}
}
   \caption{Maximum GMRES iteration counts for solving systems with matrix $\AmatoI\Amatt$, where $A_1=A_2=I$ and $\NLiDRR{n_1-n_2} \leq 0.5\times  k^{-\beta}$. Grey markers indicate GMRES reached $500$ iterations without convergence.}\label{fig:n}
\end{figure}

    \begin{figure}[h!]
    \subfigure[$\beta = 0,0.1,0.2,0.3.$]  
    {
    \scalebox{0.55}{
      \centering
\input{nbpc-linfinity-plot-A-0.pgf}
}
}
\subfigure[$\beta = 0.4,0.5 ,0.6,0.7.$]
{
        \scalebox{0.55}{
      \centering
\input{nbpc-linfinity-plot-A-1.pgf}
}
}
\subfigure[$\beta = 0.8,0.9,1.$]
{
          \scalebox{0.55}{
\input{nbpc-linfinity-plot-A-2.pgf}
}
}
\caption{Maximum GMRES iteration counts for solving systems with matrix $\AmatoI\Amatt$, where $n_1=n_2=1$ and $\NLiDRRdtd{A_1-A_2} \leq k^{-\beta}$.}\label{fig:A}
\end{figure}

    \begin{figure}[h!]
    \subfigure[ $\beta = 0,0.1,0.2,0.3.$]
{
          \scalebox{0.55}{
\input{nbpc-linfinity-plot-deterministic-10-n-0.pgf}
}
}
\subfigure[$\beta = 0.4,0.5,0.6,0.7.$]
{
          \scalebox{0.55}{
\input{nbpc-linfinity-plot-deterministic-10-n-1.pgf}
}
}
 \subfigure[$\beta = 0.8,0.9,1.$]
{
\centering{
          \scalebox{0.55}{
\input{nbpc-linfinity-plot-deterministic-10-n-2.pgf}
}
}
}
   \caption{GMRES iteration counts for solving systems with matrix $\AmatoI\Amatt$, where $A_1=A_2=I$ and $\NLiDRR{n_1-n_2} = 0.5\times  k^{-\beta}$  and $n_2$ is a $10\times10$ checkerboard. Grey markers indicate GMRES reached $500$ iterations without convergence.}\label{fig:ndet}
\end{figure}

    \begin{figure}[h!]
    \subfigure[ $\beta = 0,0.1,0.2,0.3.$]
{
          \scalebox{0.55}{
\input{nbpc-linfinity-plot-deterministic-2-n-0.pgf}
}
}
\subfigure[$\beta = 0.4,0.5,0.6,0.7.$]
{
          \scalebox{0.55}{
\input{nbpc-linfinity-plot-deterministic-2-n-1.pgf}
}
}
 \subfigure[$\beta = 0.8,0.9,1.$]
{
\centering{
          \scalebox{0.55}{
\input{nbpc-linfinity-plot-deterministic-2-n-2.pgf}
}
}
}
   \caption{GMRES iteration counts for solving systems with matrix $\AmatoI\Amatt$, where $A_1=A_2=I$ and $\NLiDRR{n_1-n_2} = 0.5\times  k^{-\beta}$ and $n_2$ is $2\times 2$ checkerboard. Grey markers indicate GMRES reached 500 iterations without convergence.}\label{fig:ndet2}
\end{figure}

\paragraph{Description of the four different $A_2$ and $n_2$ considered.}

\ben
\item $A_2=I$ and 100 realisations of randomly-chosen $n_2$ that are piecewise constant with respect to a $10\times10$ square grid on $D$ (experiments in Figure \ref{fig:n}).
\item $n_2=1$ and  100 realisations of randomly-chosen $A_2$ that are piecewise constant with respect to a $10\times10$ square grid on $D$ (experiments in Figure \ref{fig:A}).
\item $A_2=I$ and deterministic $n_2$ that are piecewise constant in a checkerboard pattern on a $10\times10$ square grid on $D$ (experiments in Figure \ref{fig:ndet}).
\item $A_2=I$ and deterministic $n_2$ that are piecewise constant in a checkerboard pattern on a $2\times 2$ square grid on $D$ (experiments in Figure \ref{fig:ndet2}).
\een

Regarding 1:
  Each instance of $n_2$ is piecewise constant with respect to a $10\times10$ square grid on $D := (0,1)^2$ with values
  chosen independently from
  the $\Unif\mleft(1-\alpha,1+\alpha\mright)$ distribution, where
  \begin{align} \label{eq:alpha}
   \alpha : = 0.5 \times k^{-\beta},  \quad \text{for} \quad  \beta \in \{0,0.1,\ldots,0.9,1\}, 
  \end{align}
  so that $\Vert n_1 - n_2 \Vert_{L^\infty(D, \mathbb{R})} \leq  0.5 \times k^{-\beta}$.

Regarding 2: On each element of the grid we set 
$$ A_2 = \left[ \begin{array}{cc}1 + a & b \\b & 1+c \end{array} \right],   $$ 
where $a,c \in  \Unif\mleft(0,\alpha\mright)$ and $\alpha $ as given in \eqref{eq:alpha}. For each draw of $a,c$ we choose $b \in \Unif(0,\delta)$, where
$\delta = \min \{ \alpha, \sqrt{(1+a)(1+c)}\} $. It follows that each draw of $A_2$ is positive definite
almost surely.  Since the elements of $A_1 - A_2 = I - A_2$ are all bounded above by $\alpha$, direct computation shows that
\beq\label{eq:A1A2}
\N{ A_1 - A_2}_{L^\infty(D; {\rm op})} \leq 2 \alpha \leq k^{-\beta}.
\eeq

Regarding 3: As in 1., $n_2$ is piecewise constant with respect to a $10\times10$ square grid on $D$, but now $n_2$ alternates between the values $1-\alpha$ and $1+\alpha$ (with $\alpha$ as in \eqref{eq:alpha}) in a checkerboard pattern. In this case, 
$\Vert n_1 - n_2 \Vert_{L^\infty(D, \mathbb{R})} =  0.5 \times k^{-\beta}$.

Regarding 4: This set up is exactly the same as in 3., except that now the square grid on $D$ is $2\times2$.
The motivation for this checkerboard is that now there are just under 16 wavelengths in each piece of the board (of width $0.5$) when $k=100$ (i.e., at the largest $k$ considered), whereas in the $10\times10$ checkerboard there is only just over one wavelength in each piece of the board at $k=100$.

In all four cases, to obtain the action of $\AmatoI$, we calculate the exact LU decomposition of $\Amato$.
We then solve the linear system \eqref{eq:pcsystem1} using standard GMRES 
(i.e.~with residual minimisation in the Euclidean norm) and zero initial guess, and record the number of GMRES iterations taken to achieve a relative residual of $10^{-5}$. In the first and second experiments, we consider  $k$ between $20$ and $100$, and in the third and fourth experiments we consider $k$ between $20$ and $200$.
        
\paragraph{Description and discussion of results.}

In all cases we see the number of GMRES iterations being bounded (over the ranges considered) for $\beta \in\{ 0.8, 0.9, 1.0\}$, and growing for $\beta \in \{0,\ldots,0.5\}$, with the growth worse for the two deterministic cases in Figures \ref{fig:ndet} and \ref{fig:ndet2}. For $\beta \in \{0.6, 0.7\}$ we see either growth or boundedness, depending on the problem.
To check if the number of samples in Figures \ref{fig:n} and \ref{fig:A} were large enough, we repeated the experiments leading to Figures \ref{fig:n} and \ref{fig:A} with  $200$ realisations of $n_2$ (as opposed to $100$ realisations in Figure \ref{fig:n}), but observed very little change in the results.

This behaviour is better than expected from the condition \eqref{eq:suffcondition}, which covers $\beta=1$. Because of the results about \Qtwo at the PDE level in Theorem \ref{thm:2} and Lemma \ref{lem:sharp}, we expect that the number of iterations will eventually grow for $\beta<1$ for $k$ sufficiently large.

\subsection{Results involving $\|A_1-A_2\|_{L^q}$ and $\|n_1-n_2\|_{L^q}$ for $2<q\leq \infty$}\label{sec:Lq}

\paragraph{Motivation.}
The results in Theorems \ref{thm:Linfty1} and \ref{thm:Linfty2} involving $\|A_1-A_2\|_{L^\infty}$ and $\|n_1-n_2\|_{L^\infty}$ give an answer to Q1 that fits well with both PDE considerations (via Theorem \ref{thm:2}) and physical intuition (via Remark \ref{rem:physical1k}), and is consistent with the numerical experiments for moderate $k$ in \S\ref{sec:num}.
However, the following example shows the disadvantage of measuring the difference in coefficients in the $L^\infty$ norm. Let  $D=\{ \bx: \|\bx\|_2 =2\}$, 
$A_1=A_2=I$, and 
\beq\label{eq:noweak}
n_1(\bx) =
\begin{dcases}
  1 &\tif \|\bx\|_2 \geq \half,\\
  1/2  &\tif \|\bx\|_2 < \half,
  \end{dcases}
\quad\tand\quad
n_2(\bx) =
\begin{dcases}
  1 &\tif \|\bx\|_2 \geq \half+\alpha,\\
  1/2  &\tif \|\bx\|_2 < \half+\alpha,
  \end{dcases}
\eeq
for some $0 < \alpha < 1/2$. Then $\NLiDRR{n_1-n_2} = 1/2$ for all $\alpha$, but one may expect that GMRES applied to $\AmatoI\Amatt$ would converge in a $k$-independent number of iterations for small $\alpha$. 
In this subsection, we therefore state analogues of Theorems \ref{thm:Linfty1} and \ref{thm:Linfty2}, with the differences in $n_1-n_2$ and $A_1-A_2$ measured in weaker norms than the $L^\infty$ norm. 

\paragraph{Statement of the results.}
To state these results, we need the definition that
\beq\label{eq:Lqnorm}
 \NLqDRRdtd{A}:= \left(\int_D \|A(\bx)\|_2^q \,\rd \bx\right)^{1/q}.
\eeq
Similar to in \S\ref{sec:Linfty}, Theorem \ref{thm:1anew} gives results for standard GMRES (i.e., using the vector Euclidean norm $\|\cdot\|_2$), whereas Theorem \ref{thm:1new} gives results in the weighted vector norms $\NDmatk{\cdot}$ and $\NDmatkI{\cdot}$.

\bth[Answer to Q1: $k$-independent weighted GMRES iterations]\label{thm:1new}

\

\noindent Assume that 
\bit
\item $\Dm$, $A_j$, and $n_j$, $j=1,2,$ satisfy Condition \ref{cond:1}, and 
\item $\Dm$, $A_j$, $n_j$, $j=1,2$, $h$ and $p$, satisfy Condition \ref{cond:2}.
\eit
Given $k_0>0$ and $q >2$,  there exist constants $\Cotilde$ and $\Cttilde$  independent of $h$ and $k$ (but dependent on $\Dm, A_1, n_1, p$, and $\kz$) such that if 
\beq\label{eq:condalt}
\Cotilde kh^{-d/q} \NLqDRRdtd{A_1-A_2} +\Cttilde  kh^{-d/q} \NLqDRR{n_1-n_2}
\leq \frac{1}{2}
\eeq
for all $k\geq k_0$, then \emph{both} weighted GMRES working in $\|\cdot\|_{\Dmat_k}$ (and the associated inner product) applied to 
\beqs
(\matrixAone)^{-1}\matrixAtwo\uvec = \fvec
\eeqs
 \emph{and} weighted GMRES working in $\|\cdot\|_{(\Dmat_k)^{-1}}$ (and the associated inner product) applied to 
\beqs
\matrixAtwo(\matrixAone)^{-1}\vvec = \fvec
\eeqs 
with initial guess the zero vector
converge in a $k$-independent number of iterations (in the sense of Definition \ref{def:GMRESconverge}) for all $k\geq k_0$.
 \enth

\bth[Answer to Q1: $k$-independent standard GMRES iterations]\label{thm:1anew}

\

\noindent Assume that 

\bit
\item $\Dm$, $A_j$, and $n_j$, $j=1,2,$ satisfy Condition \ref{cond:1}, and
\item $\Dm$, $A_j$, and $n_j$, $j=1,2$, $h$ and $p$, satisfy Condition \ref{cond:2}.
\eit
Given $k_0>0$, and $q >2$, let $\Cotilde$ and $\Cttilde$ be as in Theorem \ref{thm:1new}, and let $s_{+}>0$ and $m_{\pm}>0$ be as in Lemma \ref{lem:normequiv} (note that all these constants are independent of $k$, $h$, and $p$). Then if 
\beq\label{eq:condaalt}
\Cotilde \mleft(\frac{\splus}{\mminus}\mright) \frac{h^{-d/q}}{h} \NLqDRRdtd{A_1-A_2} + \Cttilde \mleft(\frac{\mplus}{\mminus}\mright) kh^{-d/q} \NLqDRR{n_1-n_2} \leq \half
\eeq
for all $k\geq k_0$, then standard GMRES (working in the Euclidean norm and inner product) applied to either \eqref{eq:pcsystem1} or \eqref{eq:pcsystem2} with initial guess the zero vector converges in a $k$-independent number of iterations (in the sense of Definition \ref{def:GMRESconverge}) for all $k\geq k_0$.
\enth

We make the following remarks.
\bit
\item
The constants $\Cotilde$ and $\Cttilde$ are given explicitly  in \eqref{eq:C1} and \eqref{eq:C2} in terms of  $q$ and 
\bit
\item $C^{(1)}_{\rm bound}$ given in Condition \ref{cond:1},
\item $C^{(1)}_{\rm FEM1}$ and $C^{(1)}_{\rm FEM2}$ given in Condition \ref{cond:2},
\item ${n_1}_{\max}$, ${n_1}_{\min}$, and ${A_1}_{\min}$ given in \eqref{eq:nbounds} and \eqref{eq:Abounds} (with $n$ replaced by ${n_1}$ and $A$ replaced by ${A_1}$),
\item $\Cinvs$ given in \eqref{eq:inverses}.
\eit
\item Theorems \ref{thm:Linfty1} and \ref{thm:Linfty2} can be viewed as special cases of Theorems \ref{thm:1new} and \ref{thm:1anew}, respectively, with $q$ formally set to $\infty$.
\item 
In the conditions  \eqref{eq:condalt} and \eqref{eq:condaalt} there is a trade-off between the norm used to measure $\no-\nt$ and the restriction on the magnitude of this norm. Indeed, as $q \searrow 2,$ $\no-\nt$ is measured in a weaker norm, but the condition 
on $\|n_1-n_2\|_{L^q}$ in \eqref{eq:condalt} and \eqref{eq:condaalt} becomes more restrictive since $h^{-d/q}$ becomes larger. Conversely, as $q \nearrow \infty,$ $\no-\nt$ is measured in a stronger norm, but the condition 
on $\|n_1-n_2\|_{L^q}$ in \eqref{eq:condalt} and \eqref{eq:condaalt} becomes less restrictive since $h^{-d/q}$ becomes smaller.
\item The condition \eqref{eq:cond} contains no powers of $h$, but the analogous condition  \eqref{eq:condalt} does.
The reason is that, in the proofs of Theorems \ref{thm:Linfty1} and \ref{thm:Linfty2}, we use the inequality $\|(n_1-n_2) v_h\|_{L^2(D)}\leq \|n_1-n_2\|_{L^\infty(D)}\|v_h\|_{L^2(D)}$, where $v_h\in V_{h,p}$, to ``pull out" the $L^\infty$ norm of $n_1-n_2$. In Theorems \ref{thm:1new} and \ref{thm:1anew} we want to ``pull out" the $L^q$ norm of $n_1-n_2$, so we use the inequality
$\|(n_1-n_2) v_h\|_{\LtD} \leq \|n_1-n_2\|_{L^q(D;\Rea)}\|v_h\|_{L^s(D)}$ where $1/2=1/q + 1/s$ (see \eqref{eq:keynbound2}). We then need to convert the $L^s$ norm of $v_h$ into an $L^2$ norm (essentially because the PDE \eqref{eq:pde} is posed in the $L^2$-based Sobolev space $H^1_{0,D}(D)$), and we do this using the inverse inequality \eqref{eq:inverses}, introducing powers of $h$. Similar considerations hold for $A_1-A_2$.
\item The remarks in Points (b), (c), and (d) in \S\ref{sec:Linfty} (for Theorems \ref{thm:Linfty1} and \ref{thm:Linfty2}) hold here for Theorems \ref{thm:1new} and \ref{thm:1anew}. Similarly, Lemma \ref{lem:standardvsweighted} also applies.
\item The numerical experiments in \cite[\S4.5.2]{Pe:20} consider $A_1=A_2=I$ and $n_j$ defined in a similar way to \eqref{eq:noweak} (but jumping across a half-plane) designed so that $\|n_1-n_2\|_{L^q(D)} = C k^{-\beta}$ (with $C>0$ independent of all parameters) for $\beta\in \{0, 0.1,\ldots, 0.9, 1.0\}$ for all $1\leq q<\infty$.
These experiments show similar behaviour to those in \S\ref{sec:num}, with the number of GMRES iterations levelling off over the range 
$10\leq k\leq 100$
 for $\beta\in \{0.7, 0.8, 0.9, 1.0\}$.
These experiments indicate that the conditions \eqref{eq:condalt} and \eqref{eq:condaalt} for $k$-independent GMRES iterations may not be sharp in their $q$- and $h$-dependence.
\eit

\subsection{Implications of these results for uncertainty quantification for the Helmholtz equation}\label{sec:UQ2}

The results of Theorems \ref{thm:Linfty1}, \ref{thm:Linfty2}, \ref{thm:1new}, and \ref{thm:1anew} and the numerical results in \S\ref{sec:num} show that 
$\|A_1-A_2\|$ and $\|n_1-n_2\|$ must decrease as $k$ increases for $\Amat_1$ to be a good preconditioner for $\Amat_2$.

This suggests that the idea of ``reusing" preconditioners in calculating $Q(u(\omegasl))$, $\ell=1,\ldots,N$ in \eqref{eq:samplingexp} becomes less effective as $k\tendi$.
However, numerical experiments in \cite[\S4.6.3, Page 183]{Pe:20}, with $A=I$ and $n$ given by an artificial Karhunen-Lo\`eve-type expansion, 
show that the number of samples $N$ in \eqref{eq:samplingexp} must increase with $k$ to make the error \eqref{eq:QMCerror} controllably small as $k\tendi$.
There is therefore a trade-off between how quickly $N$ increases with $k$ and how quickly the region in coefficient space for which a preconditioner is effective decreases with $k$. 

In the rest of this subsection we give some more details about the experiments in \cite[\S4.6.3]{Pe:20} about  ``reusing" preconditioners in calculating $Q(u(\omegasl))$. Before doing this, however,
we highlight that this strategy is related to \emph{mean-based preconditioning}, where a single preconditioner is calculated corresponding to the mean of the random coefficient. The initial computational work on mean-based preconditioning for the equation $-\nabla\cdot(A\nabla u)= f$ was carried out by \cite{GhKr:96},  \cite{PeGh:00}, and  \cite{Ke:04}, with theory following from \cite{PoEl:09} and \cite{ErPoSiUl:09}. In the Helmholtz context, mean-based preconditioning has been used in \cite{JiCa:09} and \cite{WaLi:19} (see the literature review in \cite[\S4.7]{Pe:20} for more discussion).

Preliminary experiments in \cite[\S4.6.3]{Pe:20} show that, at least for moderate $k$, the idea of ``reusing nearby preconditioners" can give a significant speedup  in the implementation of    Quasi-Monte-Carlo (QMC) methods for Helmholtz problems.
One of the experiments in \cite[\S4.6.3]{Pe:20} concerned the UQ of various quantities of interest for the solution of
  \eqref{eq:pde} on the unit square $D = (0,1)^2$ with an impedance boundary condition, and  with coefficients  taken to be $A=I$ and $n(\bx,y) = 1 + \sum_{j=1}^s y_j \psi_j(\bx),$ \text{with} $ \bx \in D$, $y_j  \in \text{Unif}(-1/2,1/2) \ \text{i.i.d.}$ and $ \psi_j(\bx) =  j^{-2} \cos(j  \pi x_1/4)\cos((j+1) \pi x_2/4)  $.   For the UQ,   a randomly shifted QMC rule taken from \cite{Nu} was used.

The first part of the experiment uses  extrapolation  to determine how the number of QMC
  points $N = N(k)$ need to be chosen to ensure that the error in the expected value of each of the quantities
  of interest is bounded with respect to $k$. These experiments  suggest  that for $k \in[10,60]$, $N(k)$ needs to
  grow somewhere between linearly and quadratically in $k$ to keep the error bounded; thus providing an empirical answer in this case to Q2 in \S\ref{sec:UQ1} (see \eqref{eq:QMCerror}).
  
Since $N(k)$ instances of the boundary value problem have to be solved for each $k$, \cite[\S4.6.2]{Pe:20} looks at speeding  this computation up by computing a  set of `preconditioning points'  on a
  (coarse) sublattice of the QMC points. By the theory above, for the sublattice to be guaranteed useful,  for each point $y$ on it, there should be another coarse lattice point  $y'$ satisfying approximately $\Vert n(\cdot, y ) - n(\cdot, y') \Vert_{L^\infty(D)} \lesssim k^{-1}$.  In practice, it is hard to compute such a lattice (see Q3 in \S \ref{sec:UQ1}), so in \cite[\S4.6.3]{Pe:20} a coarse sublattice of $N_c = N_c(k)$ points is
  computed that instead equidistributes the right-hand side of the upper bound  
    $$\Vert n(\cdot, y) - n(\cdot,y')\Vert_{L^\infty(D)} \ \leq \ \sum_{j=1}^s j^{-2} \vert y_j - y_j'\vert,   $$  with the value of $N_c$ again determined by extrapolation from small values of $k$. Results for this procedure are given in  \cite[Table 4.5]{Pe:20} and these show that  the ratio $N_c/N$ is never bigger than  $0.016$ for $k \in [10,60]$ in dimension $s = 10$.  Therefore, more than 98\% of the required linear systems  can be solved  by preconditioned GMRES. Moreover, the average number of GMRES iterations needed for the solution of all the required linear systems varies between 6.46 to 7.41 as $k$ increases from $10$ to $60$, showing the good quality and ``$k-$robustness'' of the set of  preconditioning points. In the experiments in \cite{Pe:20}, the restarted version of  GMRES (the default in Firedrake) was used with restarts after 30 iterations.
  This preconditioning strategy for Helmholtz UQ will be investigated in more detail  elsewhere.

\section{Definitions of the Helmholtz boundary value problem, its finite-element solution, and weighted GMRES}
\label{sec:BVP}

\paragraph{Function-space notation:} 
Let $\spd$ be the set of all symmetric-positive-definite matrices in $\RRd\times\RR^d$. 
For $D \subseteq \RRd,$ let $\LiDRRdtd$ be the set of all matrix-valued functions $A:D\rightarrow\spd$ such that $A_{i,j} \in \LiDRR$ for all $i,j =   1,\ldots,d.$ Where the range of functions is $\CC$ we suppress the second argument in a function space, e.g.~writing $\LtD$ for $\LtDCC.$

\begin{assumption}[Assumption on the domain and the coefficients]\label{ass:1}

\

(i) $\Dm$ is an open bounded Lip\-schitz set such that the open complement $\Dp \de \RRd \setminus \Dmclos$ is connected. $\Dtilde$ is a bounded connected Lipschitz open set such that $\Dmclos \compcont \Dtilde.$ Let $D \de \Dtilde \setminus \Dmclos$.

(ii) $n \in \LiDRR$ is such that $\supp\mleft(1-n\mright) \compcont \Dtilde$ and there exist $0<\nmin\leq \nmax<\infty$ such that
  \beq\label{eq:nbounds}
\nmin \leq n(\bx) \leq \nmax \quad\tfae \bx \in D,
  \eeq

(iii) $A \in \LiDRRdtd$ is such that $\supp\mleft(I-A\mright)\compcont \Dtilde$ and there exist $0<\Amin\leq \Amax<\infty$ such that
  \beq\label{eq:Abounds}
\Amin \abs{\bxi}^2 \leq \mleft(A(\bx) \bxi \mright) \cdot \bxibar \leq \Amax \abs{\bxi}^2 \quad\tfa \bxi \in \CCd \tfae \bx \in D,
\eeq
 \end{assumption} 

Let $\GD \de \partial \Dm,$ and $\GI \de \partial \Dtilde.$ (Observe that, by construction, $\partial D = \GD \cup \GI$ and $\GD \cap \GI = \emptyset.$)
Let $\gamma: \HoD\rightarrow H^{1/2}(\partial D)$ denote the Dirichlet trace operator, and $\dn:H^1(D,\Delta)\rightarrow H^{-1/2}(\partial D)$ the Neumann trace operator. 

\begin{definition}[Truncated Exterior Dirichlet Problem (TEDP)]\label{def:TEDP}
Given $\Dm$, $\Dtilde$, $A$, and $n$ satisfying Assumption \ref{ass:1}, $k>0$, $f\in \LtD$,
$\gI \in \LtGI$, and a bounded linear map $\T:\HhGI\rightarrow \HmhGI$,  we say $u \in \HoD$ satisfies the \defn{truncated exterior Dirichlet problem} if
  \beq\label{eq:tedp}
\grad \cdot \mleft(A \grad u\mright) + k^2 n u = -f \tin D,
\eeq
\beqs
\gamma u = 0, \ton \GD \tand
\eeqs
\beq\label{eq:ibc}
\dn u -T(\gamma u) = g_I \ton \GI,
\eeq
where the equation \eqref{eq:tedp} is understood in the following weak sense:
\beqs
\int_{D} 
\nabla u  \cdot(A\overline{\nabla\phi}) 
- k^2 n u\overline{\phi} 
=\int_{D} f\, \overline{\phi} \quad\tfa \phi \in C^\infty_{\comp}(D),
\eeqs
where $C^\infty_{\comp}(D):= \{ \phi \in C^\infty(D) : \supp \,\phi \text{ is a compact subset of } D\}.$
\end{definition}

Green's identity (see, e.g., \cite[Lemma 4.3]{Mc:00}) then implies that the solution of the TEDP is the solution of the following variational problem (and vice versa):
\beq\label{eq:vp}
\text{ find } u \in H_{0,D}^1(D) \text{ such that } \,a(u,v) = F(v) \,\tfa v \in H_{0,D}^1(D), 
\eeq
where
\beq\label{eq:spaceEDP}
H_{0,D}^1(D):= \big\{ v\in H^1(D) : \gamma v=0 \ton \Gamma_D\big\},
\eeq
\beq\label{eq:agen}
\aG(w,v) \de \int_{D} \mleft(\mleft(A \grad w\mright)\cdot\grad \vbar - k^2 n\minispace w \vbar\mright) - \DPGI{\T \trGI w}{\trGI v} \quad\tand \quad
F(v) = \int_D f \,\vbar + \langle g, v \rangle_{\Gamma_I},
\eeq
where $\langle\cdot,\cdot\rangle_{\Gamma_I}$ denotes the duality pairing on $\Gamma_I$ that is linear in the first argument and antilinear in the second argument.
When $A=A_j$ and $n=n_j$, for either $j=1$ or $j=2$, we denote the sesquilinear form $a(\cdot,\cdot)$ by $a_j(\cdot,\cdot)$.

\bre[Choices of the operator $T$]\label{rem:T}
We consider the two choices $T= \ri k$ and $T= \DtN$, where $\DtN : H^{1/2}(\Gamma_R) \rightarrow H^{-1/2}(\Gamma_R)$ is the Dirichlet-to-Neumann (DtN) map for the equation $\Delta u+k^2 u=0$ posed in the exterior of $\Dtilde$ with the Sommerfeld radiation condition \eqref{eq:src}; when $\Dtilde$ is a ball, the definition of $\DtN$ in terms of Hankel functions and polar coordinates (when $d=2$)/spherical polar coordinates (when $d=3$) is given in, e.g., \cite[Equations 3.5 and 3.6]{ChMo:08}, \cite[\S2.6.3]{Ne:01}. 
When $T= \DtN$, $g=0$, and $\supp f \subset\subset \Dtilde$, the solution of the TEDP equals the restriction to $D$ of the solution of the exterior Dirichlet problem posed in $D_+$ with the Sommerfeld radiation condition at infinity; see, e.g., \cite[Lemma 3.3]{GrPeSp:19}.

A simple approximation to $\DtN$ is the operator of multiplication by $\ri k$. The boundary condition \eqref{eq:ibc} then becomes the impedance boundary condition, which is also known as a first-order absorbing boundary condition (see, e.g, \cite{BaGuTu:82}, \cite[\S3.3]{Ih:98}).
See \cite[Equation 5.3]{GaMuSp:19} for analysis of the approximation $\DtN \approx \ri k$ in the case when $\Dtilde$ is a ball, and \cite{GaLaSp:21} for upper and lower bounds on the relative error between the solution of the TEDP and the exterior Dirichlet problem for more general $\Dtilde$.
\ere

\bre[Existence and uniqueness]\label{rem:existence}
With either of the choices $T= \ri k $ and $T= \DtN$, $a(\cdot,\cdot)$ satisfies a G\aa rding inequality (for $T=\ri k$ this is straightforward; for $T=\DtN$ see, e.g., \cite[Theorem 2.6.6]{Ne:01}). Therefore, once uniqueness is established the Fredholm alternative then gives existence. 
The unique continuation principle (UCP) holds (and therefore gives uniqueness) when $n\in L^\infty$ and (i) $d=2$ and $A\in L^\infty$, (ii) $d=3$ and $A\in C^{0,1}$; see the references in \cite[\S1]{GrPeSp:19}, \cite[\S2]{GrSa:20} (where the latter paper applies these results specifically to the Helmholtz equation).
The UCP does not hold in general when $d=3$ and $A$ is rougher than Lipschitz (see \cite{Fi:01}), but uniqueness is proved for a particular class of $A\in L^\infty$ in \cite[Theorem 2.7]{GrPeSp:19}.
\ere

\paragraph{Finite-element solution of the TEDP.}

Let $(\Vhp)_{h>0}$ consist of continuous piecewise-polynomials on a family of quasi-uniform simplicial meshes $\cTh$ with mesh-size $h$
and fixed polynomial degree $p$ (i.e.~we consider the ``$h$-FEM"). Note that the dimension $N$ of $\Vhp$ then satisfies $N\sim h^{-d}$, with hidden constant dependent on $p$. (The assumption of quasi-uniformity can, in principle, be relaxed, see Remark \ref{rem:ggsqu}.) 

Recall that our main results, Theorems \ref{thm:Linfty1}, \ref{thm:Linfty2}, \ref{thm:1new}, and \ref{thm:1anew}, are proved under Conditions \ref{cond:1} and \ref{cond:2}.
Some of the cases under which these conditions hold require $\partial D$ to be at least $C^{1,1}$. 
 For such $D$ it is not possible to fit $\partial D$ exactly with simplicial elements, and 
some analysis of non-conforming error is therefore necessary;  since this is standard (see, e.g., \cite[Chapter 10]{BrSc:08}), we ignore this issue here.

\begin{definition}[Finite-element approximation of the solution of the TEDP]\label{prob:fevgen}
With $a(\cdot,\cdot)$ and $F$ given in \eqref{eq:agen}, find $\uh \in \Vhp$ such that
\beq\label{eq:Galerkin}
a(\uh,\vh) = F(\vh) \tforall \vh \in \Vhp.
\eeq
We say that $\uh \in \Vhp$ is the \defn{finite-element approximation of $u$}.
\end{definition}
Conditions on $k, h$, and $p$ under which $u_h$ exists and is unique are discussed in \S \ref{sec:cond2hold}. 
(Regarding the notation $\uh$: although the finite-element approximation depends on both $h$ and $p$, since we consider decreasing $h$ with $p$ fixed in this paper, for brevity the notation only reflects the $h$-dependence.)

The matrices associated with our finite-element discretisation are defined as follows. Let $\{\phi_i, i= 1, \ldots, N\}$ be a nodal basis for $\Vhp$ with each $\phi_i$ \emph{real-valued}.
Let 
\beq\label{eq:matrixSjdef}
\big(\Smat_{A}\big)_{ij}\de \int_D \big(A \nabla \phi_j)\cdot\nabla \phi_i, \quad
\big(\Mmat_{n}\big)_{ij}\de \int_D n\,\phi_i\, \phi_j,
\quad\tand\quad
\big(\Nmat\big)_{ij}\de \int_{\GR} \T (\gamma\phi_j) \,\gamma \phi_i
\eeq
be the stiffness, domain-mass, and boundary-mass matrices, respectively. Note that both $\Smat_A$ and $\Mmat_n$ are \emph{real-valued}, but in general $\Nmat$ is \emph{complex-valued} (because both the  DtN operator $\DtN$ and the impedance operator $\ri k$ are complex-valued).
Let
\beqs
\Amat \de \Smat_{A} - k^2 \Mmat_{n} - \Nmat,
\eeqs
and let $u_h\de \sum_j u_j \phi_j$. Then \eqref{eq:Galerkin} implies that $\Amat \uvec = \fvec,$
where $(\bu)_i \de u_i$ and $(\fvec)_i \de F(\phi_i)$.
For $j=1,2$, let 
\beq\label{eq:matrixAjdef}
\Amatj \de \Smat_{A_j} - k^2 \Mmat_{n_j} - \Nmat
\eeq
and let
\beqs
\Dmat_k\de \Smat_I + k^2 \Mmat_1;
\eeqs
with $\Smat_I$ and $\Mmat_1$ the standard stiffness and mass matrices defined by \eqref{eq:matrixSjdef} above (thus $\Dmat_k$ is independent of the coefficients $A$ and $n$).
We then define the weighted norm $\|\cdot\|_{\Dmat_k}$ by \eqref{eq:Dk}. 

\paragraph{Weighted GMRES.}

We now give the set-up for weighted GMRES, introduced in \cite{Es:98} (with GMRES introduced in \cite{SaSc:86}). Consider the abstract  linear system 
$\matrixC \xvec = \dvec$
in $\mathbb{C}^N$, where $\matrixC \in \CC^{N\times N}$ is invertible. Let $\xvecz$ be the initial guess, and define the initial residual $\rvec^0 \de \matrixC \xvec^0-\dvec$ and the standard Krylov spaces:  
\beqs  
\cK^m(\Cmat, \rvec^0) \de \mathrm{span}\big\{\matrixC^j \rvec^0 : j = 0, \ldots, m-1\big\}.
\eeqs
Analogously to the definition of $\NDk{\cdot}$ above, let $(\cdot , \cdot )_{\Dmat}$ denote the inner product on $\CC^n$ 
induced by some Hermitian positive-definite matrix $\Dmat$, i.e.~$(\vvec,\wvec)_{\Dmat} \de (\Dmat \vvec, \wvec)_2,$ and let $\Vert \cdot \Vert_\Dmat$ be the induced norm. For $m \geq 1$, define the $m$th GMRES iterate $\xvec^m$  to be the unique element of $\cK^m$ such that its residual $\rvec^m:= \matrixC \xvec^m-\dvec$ satisfies  the  
 minimal residual  property: 
\beq\label{eq:minres}
 \ \Vert \rvec^m \Vert_\Dmat = \Vert  \matrixC \xvec^m-\dvec \Vert_\Dmat \ = \ \min_{\xvec \in \cK^m(\matrixC, \rvec^0)} \Vert   {\matrixC} {\xvec} -{\dvec}\Vert_\Dmat.
 \eeq
The algorithm stops when $\|\rvec^m\|_2/\|\rvec^0\|_2\leq \eps$ for some prescribed tolerance $\eps$.
Observe that when $\Dmat = \Imat$ this is the standard GMRES algorithm. We also note that in general, weighted GMRES requires the use of weighted Arnoldi process, also introduced in \cite{Es:98}, see also the alternative implementations of the weighted Arnoldi process in \cite{GuPe:14}.

%

\section{Conditions \ref{cond:1} and \ref{cond:2}}\label{sec:conditions}

\subsection{Definitions of Conditions \ref{cond:1} and \ref{cond:2}}

Condition \ref{cond:1} is a condition on the domain $D$ and coefficients $A$ and $n$ in the TEDP of Definition \ref{def:TEDP}, and involves a constant $C_{\rm bound}$ (see \eqref{eq:bound1}). When $D$, $A_j$, and ${n_j}$ satisfy Condition \ref{cond:1}, we write the corresponding constant as $C^{(j)}_{\rm bound}$.

Condition \ref{cond:2} is a condition on $h$ and $p$ in the finite-dimensional subspaces $V_{h,p}$ \emph{and} implicitly also a condition on $D$, $A$, and $n$ in the TEDP (since $D$, $A$, and $n$ must satisfy certain properties for the finite-element approximation to exist).
Condition \ref{cond:2} involves the constants $C_{\rm FEM1}$ and $C_{\rm FEM2}$; 
when $D$, $A_j$, and $n_j$ satisfy this, we write the corresponding constants as $C^{(j)}_{\rm FEM1}$ and 
$C^{(j)}_{\rm FEM2}$.

\begin{condition}[Nontrapping bound on $u$]\label{cond:1}
$D, A,$ and $n$ are such that, given $f\in L^2(D)$
  the solution of the variational problem \eqref{eq:vp} with 
  \beq\label{eq:LGf}
  F(v) = \int_D f\vbar,
  \eeq
$u$, exists, is unique, and, given $k_0>0$, $u$ satisfies the bound 
\beq\label{eq:bound1}
\big\|u\big\|_{\HokD} \leq C_{\rm bound} \N{f}_{L^2(D)} \quad \tfa k\geq k_0,
\eeq
where $C_{\rm bound}$ is independent of $k$, but dependent on $A, n, D$, and $k_0$.
\end{condition}

To state Condition \ref{cond:2} we need the definition of the norm on $(\HoDkk)^*$,
\beq\label{eq:dualnorm}
\N{F}_{(\HoDkk)^*}:= \sup_{v\in H^1_{0,D}(D)} \frac{|F(v)|}{\N{v}_{\HoDkk}}.
\eeq

\begin{condition}[$k$-independent accuracy of the FE solution for $a(\cdot,\cdot)$]
\label{cond:2}

Given $\kz>0$, 

\noindent\ben

\item $D$, $A$, $n$, $h$ and $p$ are such that, if $f= \widetilde n\sum_j \alpha_j\phi_j$ for some $\alpha_j \in \CC$ and  $\widetilde{n}\in L^\infty(D;\Rea)$  (i.e.~$f$ is an arbitrary element of $\Vhp$ multiplied by $\widetilde n$), then for all $k\geq k_0$,
  \bit
\item the solution $u_h$ of \eqref{eq:Galerkin} with $F(v)$ given by \eqref{eq:LGf} exists and is unique, and
\item the error bound
  \beq\label{eq:bound3}
\N{u-u_h}_{\HokD} \leq C_{\rm FEM1} \N{f}_{\LtD} \quad\tfa k\geq k_0, 
\eeq
holds, where $C_{\rm FEM1}$  is independent of $k$ and $h$, but dependent on $A, n, D, k_0$, and $p$.
  \eit

\item
$D$, $A$, $n$, $h$ and $p$ are such that, if $F(v)= (\widetilde{A}\nabla \widetilde{f},\nabla v)_{\LtD}$, where $\widetilde{A}\in \LiDRRdtd$ and $\widetilde{f} \de \sum_j \alpha_j \phi_j$ with $\alpha_j\in \CC$  (i.e.~$\widetilde{f}$ is an arbitrary element of $\Vhp$), then for all $k \geq \kz$,
  \bit
\item the solution $u_h$ of \eqref{eq:Galerkin} exists and is unique, and  
\item the error bound
\beq\label{eq:bound4}
\N{u-u_h}_{\HokD} \leq C^{}_{\rm FEM2}\,k\, \N{F}_{(\HoDkk)^*} \quad\tfa k\geq k_0, 
\eeq
holds, where $C^{}_{\rm FEM2}$  is independent of $k$ and $h$, but dependent on $A, n, D, k_0$, and $p$.  
  \eit
\een
\end{condition}

\subsection{When does Condition \ref{cond:1} hold?}\label{sec:cond1hold}

The summary is the following:
\ben
\item  Condition \ref{cond:1} holds when $D$, $A$, and $n$ are such that the problem is \emph{nontrapping}, in the sense that all rays starting in $D$ hit $\Gamma_I$ after some uniform time.
\item The definition of nontrapping is subtle in the case when \emph{either} $\Dm \neq\emptyset$ \emph{or} $A$ and $n$ are discontinuous.
\item Expressions for $C_{\rm bound}$ that are explicit in $A_1$ and ${n_1}$ have recently been obtained \cite{GaSpWu:20}, \cite{GrPeSp:19}, \cite{MoSp:19}; essentially $C_{\rm bound}$ is a constant multiple of the length of the longest ray in $D$.
\een

\paragraph{Regarding 1:}
In the case when $\Dm=\emptyset$, and $A$ and $n$ are $C^\infty$, the rays of the Helmholtz equation $\nabla\cdot(A\nabla u)+ k^2 nu =-f$ are defined as the projections in $\bx$ of the solutions $(\bx(s), \bxi(s)) \in \Rea^d\times \Rea^d$ of the Hamiltonian system
\beqs
\diff{x_i}{s}(s) = \pdiff{}{\xi_i}H\big(\bx(s), \bxi(s) \big), \qquad
\diff{\xi_i}{s}(s)
 = -\pdiff{}{x_i}H\big(\bx(s), \bxi(s) \big),
\eeqs
where the Hamiltonian $H(\bx,\bxi)$ given by 
\beqs
H(\bx,\bxi):= \frac{1}{n(\bx)}\sum_{i=1}^d\sum_{j=1}^{d} A_{ij}(\bx)\xi_i \xi_j - 1
\eeqs
(observe that $H$ is the \emph{semiclassical principal symbol} of the Helmholtz equation; see, e.g., \cite[\S7]{GrPeSp:19} for an introductory discussion on this).

When $T=\DtN$, the propagation of singularities results of \cite[\S VI]{DuHo:72} (see also, e.g., \cite[Chapter 24]{Ho:85}, \cite[\S12.3]{Zw:12}) combined with the parametrix argument of \cite{Va:75} then imply that Condition \ref{cond:1} is satisfied. These arguments, however, do not give an explicit expression for $C_{\rm bound}$ in terms of (properties of) $A$ and $n$.

When $T=\ri k$ the results \cite[Theorems 5.5 and 5.6 and Proposition 5.3]{BaLeRa:92} imply that Condition \ref{cond:1} is satisfied, but do not give an explicit expression for $C_{\rm bound}$ in terms of (properties of) $A$ and $n$.

\paragraph{Regarding 2:}

When $A$ and $n$ are $C^\infty$, but $\Dm \neq\emptyset$, defining how the rays interact with the boundary $\Gamma$ is subtle, and requires the notion of the Melrose--Sj{\"o}strand generalized-bicharacteristic flow \cite[\S24.3]{Ho:85}, \cite{MeSj:78}, \cite{MeSj:82}. However, under a suitable nontrapping hypothesis \cite[Definition 7.20]{MeSj:82}, the results of \cite{Va:75} then imply that Condition \ref{cond:1} is satisfied when $T=\DtN$. When $T=\ri k$, the analogous result is given in \cite{BaLeRa:92} (with the specific case of a nontrapping Dirichlet obstacle covered by \cite[Equation 5.2]{BaLeRa:92}).
When $A$ and $n$ are discontinuous, the signs of the jumps dictate whether the problem is nontrapping or trapping; see \cite{PoVo:99a}, \cite{CaPoVo:99}, and \cite{MoSp:19}.

\paragraph{Regarding 3:} 
When $T=\DtN$, $D$, $A$, and $n$ are such that the problem is nontrapping, both ${n_1}$ and $A_1$ are globally $C^{1,1}$ and $C^\infty$ in a neighbourhood of $\Dm$, and no rays are tangent to $\Gamma_D$ to infinite order then \cite[Theorem 2 and Equation 6.32]{GaSpWu:20} proves that Condition \ref{cond:1} holds with
\beqs
C_{\rm bound} := \frac{2\sqrt{2}}{\pi} \frac{1}{ (n_{\rm min})^{1/2}} L(A,n,D),
\eeqs
where $L(A,n,D)$ is the length of the longest ray in a fixed neighbourhood of $D$;
this result also holds when $\Dm=\emptyset$ and both ${n_1}$ and $A_1$ are globally $C^{1,1}$.
(When $A_1=I$ and $n_1=1$ the rays are straight lines, and so no rays being tangent to $\Gamma_D$ to infinite order means that the boundary is never flat to infinite order.)

When (i) $T=\DtN$ and $\Dm$ is star-shaped with respect to a point, or (ii) $T=\ri k$ and $\Dm$ and $\Dtilde$ are both star-shaped with respect to the same point, then expressions for $C_{\rm bound}$ are given in \cite{GrPeSp:19} when $A$ and $n$ satisfy certain monotonicity conditions in the radial direction. These monotonicity conditions allow for $A$ and $n$ to be discontinuous (with nontrapping jumps), and the results of 
 \cite{GrPeSp:19} thus recover the earlier results of \cite{MoSp:19} for piecewise constant $A$ and $n$ when $T=\DtN$.

\subsection{When does Condition \ref{cond:2} hold?}\label{sec:cond2hold}

\paragraph{Part 1 of Condition \ref{cond:2}:}

When $D$, $A$, and $n$ satisfy both Condition \ref{cond:1} and additional smoothness requirements, Part 1 of Condition \ref{cond:2} holds (at least when $T= \ri k$) when $k(hk)^{2p}\leq \cC_1$ for some $\cC_1$ depending on $C_{\rm bound}$.

Indeed, by \cite[Theorem 2.39]{Pe:20}, when $D$, $A$, and $n$ satisfy Condition \ref{cond:1}, $T=\ri k$, 
$\partial D\in C^{p,1}$, $A \in C^{p-1,1}$, and $n \in H^{\max\{{p-1}, \ceil{d/2} + 1\}}$ 
there exist $\cC_2, \cC_3>0$, depending on $A$ and $n$ (but not on $C_{\rm bound}$), such that if 
\beq\label{eq:threshold1}
k(hk)^{2p}\leq  \cC_2 (C_{\rm bound})^{-1},
\eeq
then Part 1 of Condition \ref{cond:1} holds with $C_{\rm FEM1} =\cC_3 C_{\rm bound}$.
(We remark that the condition $n \in H^{\max\{{p-1}, \ceil{d/2} + 1\}}$ is imposed in \cite[Theorem 2.39]{Pe:20} so that 
if $v \in H^{p-1}(D)$, then $n v \in H^{p-1}(D)$ \cite[Page 63, paragraph under Assumption 2.35]{Pe:20}; by \cite[Theorem 1.4.1.1, Page 21]{Gr:85} this property is also ensured by requiring $n\in C^{\min\{p-2,0\},1}(D)$.)

Furthermore, the arguments in \cite{LaSpWu:19a} show that, when $D$, $A$, and $n$ satisfy Condition \ref{cond:1}, $T=\DtN$, $p=1$, $\partial D\in C^{0,1}$, $A\in C^{0,1}$, and $n\in L^\infty$,  
there exist $\cC_2, \cC_3>0$, depending on $A$ and $n$ (but not on $C_{\rm bound}$), such that if \eqref{eq:threshold1} holds, 
then Part 1 of Condition \ref{cond:1} holds with $C_{\rm FEM1} =\cC_3 C_{\rm bound}$. 
Indeed, while \cite[Theorem 1.5]{LaSpWu:19a} obtains a bound on the \emph{relative error} of the finite-element approximation (rather than the error in terms of the data as in \eqref{eq:bound3}), using the bound $\|u\|_{H^2}\lesssim k \, C_{\rm bound}\|f\|_{L^2}$ in \cite[Equation 2.5]{LaSpWu:19a} gives \eqref{eq:bound3}, 
with this bound on the $H^2$ norm following from \eqref{eq:bound1} by elliptic regularity; see, e.g., \cite[Remark 2.14]{GrPeSp:19} and the references therein.

The arguments in both \cite{Pe:20} and \cite{LaSpWu:19a} use the ``elliptic-projection" modification of the Schatz duality argument from \cite{Wu:14}, \cite{ZhWu:13}, \cite{DuWu:15}, with these latter papers pioneering this argument for the case $A=I$ and $n=1$ (and \cite{DuWu:15}, and subsequently \cite{Pe:20}, using an additional ``error-splitting" technique).

\paragraph{Part 2 of Condition \ref{cond:2}:}

When $D$, $A$, and $n$ satisfy both Condition \ref{cond:1} and additional smoothness requirements, Part 2 of Condition \ref{cond:2} holds (at least when $T= \ri k$) when $k(hk)^p \leq \cC_4$ for some $\cC_4$ depending on $C_{\rm bound}$.

Indeed, by \cite[Theorem 2.15 and Remark 2.2]{ChNi:19},
when $D$, $A$, and $n$ satisfy Condition \ref{cond:1}, $T=\ri k$, $\partial D\in C^{p,1}$, $A \in C^{p-1,1}$, and $n \in C^{p-1,1}$
then there exists $\cC_4, \cC_5>0$, depending on $A$, $n$ (but not $C_{\rm bound}$) such that if 
\beqs
k(hk)^p \leq \cC_4 (C_{\rm bound})^{-1},
\eeqs
then 
\beq\label{eq:qo}
\N{u-u_h}_{H^1_k(D)}\leq \cC_5 \min_{v_h\in V_{h,p}}\N{u-v_h}_{H^1_k(D)}
\eeq
i.e.~the finite-element approximation is \emph{quasi-optimal}.
Lemma \ref{lem:H1} shows that, if Condition \ref{cond:1} holds, then $\|u\|_{\HoDkk}\lesssim k \|F\|_{(\HoDkk)^*}$, and then, setting $v_h=0$ in \eqref{eq:qo}, we see that \eqref{eq:bound4} holds.

\section{Proofs of Theorems \ref{thm:Linfty1}, \ref{thm:Linfty2}, \ref{thm:1new} and \ref{thm:1anew}
and Lemma \ref{lem:standardvsweighted}}\label{sec:proofs}

\subsection{Preliminary results about the variational formulation and finite-element method}

\paragraph{Inverse inequality.} Recall that for $s\geq 2$ there exists $\Cinvs>0$ (depending on the constant in the definition of quasi-uniformity -- see \cite[Equation 4.4.15]{BrSc:08}) such that
\beq\label{eq:inverses}
\NLsD{v_{h,p}} \leq \Cinvs h^{d\mleft(\frac1{s} - \half\mright)} \NLtD{v_{h,p}}\quad\tfa v_{h,p} \in V_{h,p};
\eeq
see \cite[Theorem 4.5.11 and Remark 4.5.20]{BrSc:08}; observe that if $s=2$ then $\Cinvs =1$.

\ble[Norm equivalences of FE functions]\label{lem:normequiv}
There exist $m_\pm >0$ and $s_+>0$, independent of $h$ (but dependent on $p$), such that
\beq\label{eq:normequiv1}
m_- h^{d/2} \N{\vvec}_2 \leq \N{v_h}_{\LtD} \leq m_+ h^{d/2} \N{\vvec}_2,
\eeq
and
\beq\label{eq:normequiv2}
\N{\nabla v_h}_{\LtD} \leq s_+ h^{d/2-1} \N{\vvec}_2,
\eeq
for all finite-element functions $v_h =\sum_i v_i \phi_i \in \Vhp$.
\ele

\bpf[Sketch proof]
The inequalities \eqref{eq:normequiv1} can be proved by direct computation, and then \eqref{eq:normequiv2} obtained from \eqref{eq:normequiv1} by the  inverse inequality 
\beq\label{eq:inverse2}
\|\nabla v_h\|_{L^2}\lesssim h^{-1}\|v_h\|_{L^2}
\eeq
\cite[Theorem 4.5.11 and Remark 4.5.20]{BrSc:08}. For the details, see \cite[Proof of Lemma 4.6]{Pe:20}. 
\epf

\bre[Relaxing the assumption of quasi-uniformity]\label{rem:ggsqu}
We assume that $(\Th)_{h>0}$ is a quasi-uniform family of meshes so that
we can use the inverse inequalities \eqref{eq:inverses} and \eqref{eq:inverse2}. Nevertheless,

(i)  We see below that the proof of Theorem \ref{thm:Linfty1} does not use either \eqref{eq:inverses} or \eqref{eq:inverse2}, and so Theorem \ref{thm:Linfty1} holds without the quasi-uniformity assumption.

(ii) We expect that analogues of some of the results in Theorems \ref{thm:Linfty2}, \ref{thm:1new}, and \ref{thm:1anew} can be obtained for shape-regular meshes, following the techniques in the proofs \cite[Section 3.4 and 4.1.2]{GaGrSp:15} (see Remark \ref{rem:GaGrSp} for discussion on how the results of \cite{GaGrSp:15} are related to those in the present paper). Indeed, \cite{GaGrSp:15} contains bounds on preconditioned mass matrices, analogous to those in Lemma \ref{lem:keylemma1}. From the way Lemma \ref{lem:keylemma1} is used in the proofs of Theorems \ref{thm:Linfty2}, \ref{thm:1new}, and \ref{thm:1anew}, we  expect that analogues of these results with $A_1=A_2$ hold in the case of shape-regular meshes. However, \cite{GaGrSp:15} does not contain any results on preconditioned stiffness matrices analogous to Lemma \ref{lem:keylemma2}. Therefore,
from the way Lemma \ref{lem:keylemma2} is used in the proofs of Theorems \ref{thm:Linfty2}, \ref{thm:1new}, and \ref{thm:1anew}, 
 it is unclear at this point whether analogues of these results with $A_1\neq A_2$ can be obtained in the case of shape-regular meshes.
\ere
%

\subsection{Proofs of Theorems \ref{thm:1new} and \ref{thm:1anew}}\label{sec:proofsGMRES}

In this subsection we prove Theorems \ref{thm:1new} and \ref{thm:1anew}; in the next subsection we indicate how the proofs can be modified slightly to prove Theorems \ref{thm:Linfty1} and \ref{thm:Linfty2}.

The main ingredient of the proofs of  Theorems \ref{thm:1new} and \ref{thm:1anew} is the following lemma, which we prove in \S\ref{sec:proofmi}.

\begin{lemma}\label{lem:mainingredient}
\noindent Assume that 
\bit
\item $\Dm$, ${A_1}$, and ${n_1}$, satisfy Condition \ref{cond:1}, and 
\item $\Dm$, ${A_1}$, ${n_1}$, $h$ and $p$, satisfy Condition \ref{cond:2}.
\eit
Let $\mpm$ and $\splus$ be given as in Lemma \ref{lem:normequiv}.
Given $\kz>0$ and $q >2$, let
\beq\label{eq:C1tilde}
\Cotilde \de \Cinvs\Co \quad \tand\quad \Cttilde \de \Cinvs C_2
\eeq
where 
\beq\label{eq:C1}
C_1\de
\left[ C_{\rm FEM2}^{(1)} + 
 \frac{1}{\min\big\{{A_1}_{\min},{n_1}_{\min}\big\}}\left( \frac{1}{k_0} + 2 C^{(1)}_{\rm bound}{n_1}_{\max}  \right) \right],
\eeq
\beq\label{eq:C2}
C_2:=\big(C_{\rm FEM1}^{(1)} + C_{\rm bound}^{(1)}\big),
\eeq
and $\Cinvs$ is defined by \eqref{eq:inverses} with $1/s = 1/2 - 1/q.$

Then, for all $k\geq k_0$,
\begin{align}\nonumber
&\max\set{\NDk{\Imat - \AmatoI\Amatt},\NDkI{\Imat -\Amatt\AmatoI}}\\
&\hspace{3cm} 
\leq \Cotilde kh^{-d/q} \NLqDRRdtd{A_1-A_2} + \Cttilde  kh^{-d/q}  \NLqDRR{n_1-n_2}
\label{eq:main1alt}
\end{align}
and 
\begin{align}\nonumber
&\max\set{\Nt{\Imat - \AmatoI\Amatt}, \Nt{\Imat -\Amatt\AmatoI}}\\
&\hspace{1cm}
\leq \Cotilde \mleft(\frac{\splus}{\mminus}\mright) h^{-d/q-1}\NLqDRRdtd{A_1-A_2} + \Cttilde \mleft(\frac{\mplus}{\mminus}\mright) kh^{-d/q}\NLqDRR{n_1-n_2}.
\label{eq:main1aalt}
\end{align}
\end{lemma}

\subsubsection{Proofs of Theorems \ref{thm:1new} and \ref{thm:1anew} assuming Lemma \ref{lem:mainingredient}}\label{sec:mainproofs}

Let 
\beqs
W_\Dmat(\matrixC)\de \Big\{ (\matrixC \xvec, \xvec)_{\Dmat} : \xvec \in \CCN, \|\xvec\|_\Dmat=1\Big\};
\eeqs
$W_\Dmat(\matrixC)$ is called the \emph{numerical range} or \emph{field of values} of $\matrixC$ (in the $(\cdot,\cdot)_\Dmat$ inner product).

\begin{theorem}[Elman estimate for weighted GMRES]\label{thm:GMRES1_intro} 
Let $\matrixC$ be a matrix with $0\notin W_\Dmat(\matrixC)$. Let $\beta\in[0,\pi/2)$ be defined such that
\beq\label{eq:cosbeta}
\cos \beta \de \frac{\mathrm{dist}\big(0, W_\Dmat(\matrixC)\big)}{\N{\matrixC}_{\Dmat}}.
\eeq
If the matrix equation $\matrixC \xvec = \by$ is solved using weighted GMRES then, 
for $m\in \mathbb{N}$, the GMRES residual $\rvecm$ 
satisfies
\beq\label{eq:Elman}
\frac{\N{\rvecm}_{\Dmat}}{\N{\rvecz}_{\Dmat}} \leq \sin^m \beta. 
\eeq
\end{theorem}
The bound \eqref{eq:Elman} with $\Dmat=\Imat$ was first proved in \cite[Theorem 6.3]{El:82} (see also \cite[Theorem 3.3]{EiElSc:83}) and was written in the above form in \cite[Equation 1.2]{BeGoTy:06}. The bound \eqref{eq:Elman} for arbitrary Hermitian positive-definite $\Dmat$ was stated implicitly (without proof) in \cite[p. 247]{CaWi:92} and proved in \cite[Theorem 5.1]{GrSpVa:17}.

Theorem \ref{thm:GMRES1_intro} has the following corollary, and the proofs of Theorems \ref{thm:1new} and \ref{thm:1anew}
 follow from combining this with Lemma \ref{lem:mainingredient}.

\begin{corollary}
\label{cor:GMRES_intro} 
If $\|\Imat - \matrixC \|_\Dmat \leq \alpha < 1$, then, with $\beta$ defined as in \eqref{eq:cosbeta},
\beqs
\cos \beta \geq \frac{1-\alpha}{1+\alpha}\quad\tand\quad
\sin \beta \leq \frac{2 \sqrt{\alpha}}{(1+\alpha)}.
\eeqs
\end{corollary}

\bpf[Proof of Theorem \ref{thm:1new}]
We first prove the bound for the left-preconditioned system \eqref{eq:pcsystem1} (i.e.~involving $(\Amato)^{-1} \Amatt$). 
By \eqref{eq:condalt} and Lemma \ref{lem:mainingredient} we can apply Corollary \ref{cor:GMRES_intro} with $\matrixC= (\Amato)^{-1} \Amatt$, $\Dmat=\Dmat_k$, and $\alpha=1/2$. 
Therefore 
\beq\label{eq:GMREScm}
\big\| \rvec^m\big\|_{\Dmat_k} \leq c^m \big\|\rvec^0\big\|_{\Dmat_k},
\eeq
where $c:= 2\sqrt{2}/3<1$. Therefore, given $\eps>0$, if $m\geq \log(1/\eps)/ \log(1/c)$, then $\| \rvec^m\|_{\Dmat} \leq \eps \|\rvec^0\|_{\Dmat}$.
By definition of $\rvec^m$ and the fact that $\bx^0=\bze$, 
\beqs
\rvec^m = (\Amato)^{-1} \Amatt ( \uvec^m- \uvec) \quad\tand\quad \rvec^0 = -(\Amato)^{-1} \Amatt  \uvec.
\eeqs
Therefore, 
\begin{align}\nonumber
\N{\uvec^m - \uvec}_{\Dmat_k} \leq \N{  (\Amatt)^{-1} \Amato }_{\Dmat_k} \N{\rvec^m}_{\Dmat_k}
& \leq  \N{  (\Amatt)^{-1} \Amato }_{\Dmat_k} \eps\N{\rvec^0}_{\Dmat_k},\\
& \leq  \N{  (\Amatt)^{-1} \Amato }_{\Dmat_k} \eps\N{  (\Amato)^{-1} \Amatt }_{\Dmat_k} \N{\uvec}_{\Dmat_k}.\label{eq:interim1}
\end{align}
The bound \eqref{eq:main1alt} in Lemma \ref{lem:mainingredient} implies that $\|  (\Amato)^{-1} \Amatt \|_{\Dmat_k} \leq 1 + \alpha = 3/2$. 
To obtain a bound on $\|  (\Amatt)^{-1} \Amato \|_{\Dmat_k}$, we observe that,
since $D$, $A_2$, ${n_2}$, $h$, and $p$ are assumed to satisfy Conditions \ref{cond:1} and \ref{cond:2}, the argument of Lemma \ref{lem:mainingredient} can be repeated with the indices $1$ and $2$ swapped, resulting in the bound $\|  (\Amatt)^{-1} \Amato \|_{\Dmat_k} \leq 3/2$. Inputting these bounds on $\|  (\Amato)^{-1} \Amatt \|_{\Dmat_k}$ and $\|  (\Amatt)^{-1} \Amato \|_{\Dmat_k}$ into \eqref{eq:interim1}, the result follows.

The bound for the right-preconditioned system  \eqref{eq:pcsystem2} follows in an analogous way, with $\matrixC= \Amatt(\Amato)^{-1} $, $\Dmat=(\Dmat_k)^{-1}$, and $\alpha=1/2$.
\epf

\

\bpf[Proof of Theorem \ref{thm:1anew}]
This follows from Lemma \ref{lem:mainingredient} in an analogous way to the proof of Theorem \ref{thm:1new}, except with $\Dmat=\Imat$ and using \eqref{eq:main1aalt} instead of \eqref{eq:main1alt}.
\epf

\bre[The improvement of the Elman estimate \eqref{eq:Elman} in \cite{BeGoTy:06}]
A stronger result than \eqref{eq:Elman} is given for standard (unweighted) GMRES in \cite[Theorem 2.1]{BeGoTy:06}, and then converted to a result about weighted GMRES in \cite[Theorem 5.3]{BoDoGrSpTo:19}. In the bound of  \cite[Theorem 2.1]{BeGoTy:06}, the convergence factor $\sin \beta$ is replaced by a function of $\beta$ strictly less than $\sin\beta$ for $\beta\in (0,\pi/2)$. Using this stronger result, however, does not improve the $k$-dependence in the bounds of Theorems \ref{thm:1new} or \ref{thm:1anew}.
\ere

\subsubsection{Proof of Lemma \ref{lem:mainingredient}}\label{sec:proofmi}

To prove Lemma \ref{lem:mainingredient} we need the following two lemmas.

\ble[Bounds on $(\Amato)^{-1} \Mmat_{n}$ and $\Mmatn\AmatoI$]\label{lem:keylemma1}
Assume that $D, A_1,$ and $n_1$ satisfy Condition \ref{cond:1} and that $D, A_1, n_1,h,$ and $p$ satisfy Part 1 of Condition \ref{cond:2}. 
Let $\mpm$ be as in Lemma \ref{lem:normequiv} and let $\Cotilde$ and $\Cttilde$ be given by \eqref{eq:C1tilde}.
Then, for all $n\in \LiDRR$,
$q > 2$, and $k\geq \kz$,
\beq\label{eq:keybound12}
\max\set{\NDk{\AmatoI \Mmatn},\NDkI{\Mmatn\AmatoI}} \leq \Cttilde h^{-d/q} \frac{\NLqDRR{n}}k
\eeq
and 
\beq\label{eq:keybound1a2}
\max\set{\Nt{\AmatoI \Mmatn},\Nt{\Mmatn\AmatoI}} \leq \Cttilde\mleft(\frac{\mplus}{\mminus}\mright) h^{-d/q} \frac{\NLqDRR{n}}k.
\eeq

\ele

\ble[Bounds on $(\Amato)^{-1} \Smat_A$ and $\SmatA\AmatoI$]\label{lem:keylemma2}
Assume that 
$D, A_1,$ and $n_1$ satisfy
Condition \ref{cond:1} and that
$D, A_1, n_1,h,$ and $p$ satisfy
Part 2 of Condition \ref{cond:2}. 
Let $\mpm$ and $s_+$ be as in Lemma \ref{lem:normequiv} and let $\Cotilde$ and $\Cttilde$ be given by \eqref{eq:C1tilde}.
Then for all $A\in L^\infty(D,\spd)$, $q > 2$, and $k\geq \kz$,
\beq\label{eq:keybound22}
\max\set{\NDk{\AmatoI \SmatA},\NDkI{\SmatA\AmatoI}} \leq \Cotilde h^{-d/q}k \NLqDRRdtd{A}
\eeq
and
\beq\label{eq:keybound2a2}
\max\set{\Nt{\AmatoI \SmatA},\Nt{\SmatA\AmatoI}} \leq \Cotilde\mleft(\frac{\splus}{\mminus}\mright) h^{-d/q-1} \NLqDRR{A}.
\eeq
\ele

\bpf[Proof of Lemma \ref{lem:mainingredient} using Lemmas \ref{lem:keylemma1} and \ref{lem:keylemma2}]
Using the definition of the matrices $\Amatj, \SmatA$, and $\Mmatn$ in \eqref{eq:matrixAjdef} and \eqref{eq:matrixSjdef}, we have
\begin{align}\nonumber
\Imat - (\Amato)^{-1}\Amatt = (\Amato)^{-1}\big(\Amato-\Amatt\big) &=  (\Amato)^{-1}\left( \Smat_{{A_1}} - \Smat_{{A_2}} - k^2 \big(\Mmat_{{n_1}}-\Mmat_{{n_2}}\big)\right)\\
&= (\Amato)^{-1}\left( \Smat_{{A_1}-{A_2}} - k^2 \Mmat_{{n_1}-{n_2}}\right),\label{eq:idea1}
\end{align}
and similarly 
\beq\label{eq:idea2}
\Imat -\Amatt  (\Amato)^{-1}= \left( \Smat_{{A_1}-{A_2}} - k^2 \Mmat_{{n_1}-{n_2}}\right)(\Amato)^{-1}.
\eeq
The bounds in  \eqref{eq:main1alt} on $\NDk{\Imat - (\Amato)^{-1}\Amatt}$ and  $\NDkI{\Imat - \Amatt(\Amato)^{-1}}$ then follow from using the bounds \eqref{eq:keybound12}, \eqref{eq:keybound22} in \eqref{eq:idea1} and \eqref{eq:idea2}. The bounds in \eqref{eq:main1aalt} on $\Nt{\Imat - (\Amato)^{-1}\Amatt}$ and  $\Nt{\Imat - \Amatt(\Amato)^{-1}}$ follow completely analagously, except we use the bounds \eqref{eq:keybound1a2} and \eqref{eq:keybound2a2}.
\epf

\bre[Bounds for right-preconditioning]
For brevity, we only prove the bounds on $(\Amato)^{-1} \Mmat_{n}$ in Lemma \ref{lem:keylemma1} and not those on $\Mmatn\AmatoI$. Similarly, 
we only prove the bounds on $\AmatoI\SmatA$ in Lemma \ref{lem:keylemma2} and not those on $\SmatA\AmatoI$. 
The proofs of the bounds on $\Mmatn\AmatoI$ can be found in \cite[Proof of Lemma 4.19]{Pe:20}, and the proofs of the bounds on $\SmatA\AmatoI$ can be found in \cite[Proof of Lemma 4.20]{Pe:20}. 

These proofs rely on the fact that, from the definitions of $\|\cdot\|_{\Dmat_k}$ and $\|\cdot\|_{(\Dmat_k)^{-1}}$ in \eqref{eq:Dk}, for any matrix $\Cmat \in \CCNtN$, $\|\matrixC\|_{\Dmat_k}= \|\matrixC^\dagger\|_{(\Dmat_k)^{-1}}$, where $\matrixC^\dagger$ is the conjugate transpose of $\matrixC$ (i.e.~the adjoint with respect to $(\cdot,\cdot)_2$).
This result follows from the fact that, for all $\vvec_j\in \CC^N$,
\beqs
\frac{
(\matrixC \bw_1, \bw_2)_{\Dmat_k}
}{
\N{\bw_1}_{\Dmat_k}\N{\bw_2}_{\Dmat_k}
} = 
\frac{
(\bv_1, \matrixC^\dagger \bv_2)_{(\Dmat_k)^{-1}}
}{
\N{\bv_1}_{(\Dmat_k)^{-1}}\N{\bv_2}_{(\Dmat_k)^{-1}}
} 
\eeqs
where $\wvec_j \de (\Dmat_k)^{-1}\vvec_j$, $j=1,2,$.
Therefore, since $\Mmat_n$ and $\Smat_A$ are real, symmetric matrices
\beqs
 \big\|\Mmat_n (\Amato)^{-1}\big\|_{(\Dmat_k)^{-1}}=\big\|((\Amato)^\dagger)^{-1}\Mmat_n\big\|_{\Dmat_k}
\quad\tand\quad
 \big\|\Smat_A (\Amato)^{-1}\big\|_{(\Dmat_k)^{-1}}=\big\|\big((\Amato)^\dagger\big)^{-1}\Smat_A\big\|_{\Dmat_k}.
 \eeqs 
 The matrix $(\Amato)^\dagger$ is the Galerkin matrix corresponding to the adjoint variational problem to \eqref{eq:vp}. The key point is that if Conditions \ref{cond:1} and \ref{cond:2} are satisfied for the variational problem \eqref{eq:vp}, then one can also show that they are satisfied for the adjoint problem. Therefore, the bound in \eqref{eq:keybound12} on  $\|\Mmat_n(\Amato)^{-1}\|_{(\Dmat_k)^{-1}}$ follows from the bound on $\|(\Amato)^{-1}\Mmat_n\|_{\Dmat_k}$, and the bound in \eqref{eq:keybound22} on $\|(\Amato)^{-1}\Smat_A\|_{\Dmat_k}$ also holds for $\|((\Amato)^\dagger)^{-1}\Smat_A\|_{\Dmat_k}$.
 This type of duality argument was first used in \cite{GaGrSp:15} and \cite{GrSpVa:17}.
\ere

\bpf[Proof of the bounds on $(\Amato)^{-1} \Mmat_{n}$ in Lemma \ref{lem:keylemma1}]
We first concentrate on proving \eqref{eq:keybound12}.
Given $\fvec \in \CC^N$ and $n\in \LiDRR$, we create a variational problem whose Galerkin discretisation leads to the equation $\Amato \tbu = \Mmat_n\,\fvec$.
Indeed, let $\widetilde{f} \de \sum_j f_j \phi_j\in \HozDD$. Define $\widetilde{u}$ to be the solution of the variational problem 
\beq\label{eq:411}
a_1(\widetilde{u},v)= \big(n\widetilde{f},v)_{L^2(D)} \quad\text{ for all } v\in \HozDD,
\eeq
and let $\tu_h$ be the solution of the finite-element approximation of \eqref{eq:411}, i.e.,
\beq\label{eq:41}
a_1(\tu_h,v_h)=\big(n\widetilde{f},v_h)_{L^2(D)}\quad\text{ for all } v_h\in \Vhp,
\eeq
and let $\tbu$ be the vector of nodal values of $\tu_h$. The definition of $\widetilde{f}$ then implies that \eqref{eq:41} is equivalent to the linear system $\Amato \tbu = \Mmat_{n}\,\fvec$. Since $\fvec \in \CCN$ was arbitrary, to obtain a bound on $\|(\Amato)^{-1}\Mmat_n\|_{\Dmat_k}$ we need to bound $\|\tbu\|_{\Dmat_k}$ in terms of $\|\fvec\|_{\Dmat_k}$. Because of the definition \eqref{eq:Dk}
of $\|\cdot\|_{\Dmat_k}$, this is equivalent to bounding $\|\tu_h\|_{\HokD}$ in terms of $\|\widetilde{f}\|_{\HokD}$.

Using the triangle inequality and the bounds \eqref{eq:bound3} and \eqref{eq:bound1} from Conditions \ref{cond:2} and \ref{cond:1} respectively, we find
\begin{align}
\N{\tu_h}_{\HokD} \leq
\N{\tu-\tu_h}_{\HokD} + \N{\tu}_{\HokD} \label{eq:mainevent1}
& \leq \big(C^{(1)}_{\rm FEM1}
+ C^{(1)}_{\rm bound}\big)\big\|n\ftilde\big\|_{\LtD}
\end{align}
H\"older's inequality implies that, if $q,s > 2$ such that $1/2 = 1/q+1/s,$ then
\beq\label{eq:keynbound2}
\big\|n\ftilde\big\|_{\LtD} \leq \NLqDRR{n}\big\|\ftilde\big\|_{L^s(D)}.
\eeq
Since $\ftilde \in \Vhp$, we can use the inverse inequality \eqref{eq:inverses} in \eqref{eq:keynbound2}, and combining the result with \eqref{eq:mainevent1} we obtain that
\begin{align}\label{eq:mainevent1Euan}
\N{\tu_h}_{\HokD} 
&\leq \big(C^{(1)}_{\rm FEM1}
+ C^{(1)}_{\rm bound}\big)\Cinvs \NLqDRR{n} h^{d\mleft(\frac1{s} - \half\mright)} \big\|\ftilde\|_{\LtD},\\
&\leq \big(C^{(1)}_{\rm FEM1}
+ C^{(1)}_{\rm bound}\big)\Cinvs \NLqDRR{n} h^{d\mleft(\frac1{s} - \half\mright)} k^{-1}\big\|\ftilde\big\|_{\HokD}.
\nonumber
\end{align}
The bound on $\|(\Amato)^{-1}\Mmat_n\|_{\Dmat_k}$ in \eqref{eq:keybound12} then follows from the definition \eqref{eq:Dk} of $\|\cdot\|_{\Dmat_k}$   and the definition \eqref{eq:C1tilde} of $\Cttilde$.

To prove the bound on  $\|(\Amato)^{-1}\Mmat_n\|_{2}$ in \eqref{eq:keybound1a2}, we use the consequences of \eqref{eq:normequiv1}
\beqs
m_- h^{d/2} k \N{\tbu}_2 \leq k \N{\widetilde{u}_h}_{\LtD} \leq \N{\widetilde{u}_h}_{\HokD}
\,\tand\,
\big\|\widetilde{f}\big\|_{\LtD} \leq m_+ h^{d/2}\N{\fvec}_2,
\eeqs
on either side of the inequality \eqref{eq:mainevent1Euan}. 
\epf

\

The proof of the bounds on $\AmatoI\SmatA$ in Lemma \ref{lem:keylemma2} uses the following lemma, which one can prove using fact that the sesquilinear form $a$ satisfies a  G\aa rding inequality; see \cite[Lemma 5.1]{GrPeSp:19}.

\ble[Bound for data in $\HozDDs$]\label{lem:H1}
Given $\Ftilde\in \HozDDs$, let $\widetilde{u}$ be the solution of the variational problem
\beqs
\text{ find } \,\,\widetilde{u} \in H^1_{0,D}(D) \,\,\tst \,\,
a_1(\widetilde{u},v)=\Ftilde(v) \,\, \tfa v\in H^1_{0,D}(D).
\eeqs
If Condition \ref{cond:1} holds, then $\widetilde{u}$ exists, is unique, and satisfies the bound
\beq\label{eq:bound2}
\N{\widetilde{u}}_{\HokD} \leq \frac{1}{\min\{{A_1}_{\min},{n_1}_{\min}\}}\left( 1 + 2 C^{(1)}_{\rm bound}{n_1}_{\max}  k\right) \big\|\Ftilde\big\|_{(\HokD)^*}
\eeq
for all $k\geq k_0$.
\ele

\bpf[Proof of the bounds on $\AmatoI\SmatA$ in Lemma \ref{lem:keylemma2}]
In a similar way to the proof of Lemma \ref{lem:keylemma1}, given $\fvec \in \CC^N$ and $A\in \LiDRRdtd$, let $\widetilde{f} \de \sum_j f_j \phi_j$ and observe that $\widetilde{f} \in \HozDD$. Define $\widetilde{u}$ to be the solution of the variational problem 
\beq\label{eq:411a}
a_1(\widetilde{u},v)= \Ftilde(v) \quad\text{ for all } v\in \HozDD,
\quad\text{ where } \quad
 \Ftilde(v) \de \big(A\nabla\widetilde{f},\nabla v\big)_{L^2(D)}.
\eeq
Observe that the definitions \eqref{eq:dualnorm} and \eqref{eq:1knorm} of the norms $\|\cdot\|_{(\HokD)^*}$ and $\|\cdot\|_{\HokD}$, respectively,  and the Cauchy-Schwarz inequality imply that
\beqs
\big\|\Ftilde\big\|_{(\HokD)^*}\leq \big\|A\nabla \widetilde{f}\big\|_{\LtD}.
\eeqs
With $\NLqDRRdtd{\cdot}$ defined by \eqref{eq:Lqnorm}, H\"older's inequality implies that, with 
$q,s > 2$ such that $1/2 = 1/q+1/s,$
\beqs
\big\|A\nabla \widetilde{f}\big\|_{L^2(D)}\leq \NLqDRRdtd{A} \big\|\nabla \widetilde{f}\big\|_{L^s(D)}.
\eeqs
Combining this with the inverse inequality \eqref{eq:inverses} we then have that
\beq\label{eq:keyAboundfinal}
\big\|\Ftilde\big\|_{(\HokD)^*} \leq \NLqDRRdtd{A}\big\|\grad\ftilde\big\|_{L^s(D)}\leq\Cinvs \NLqDRRdtd{A} h^{d\mleft(\frac1{s} - \half\mright)} \big\|\grad\ftilde\big\|_{L^2(D)}.
\eeq
Let $\tu_h$ be the solution of the finite-element approximation of \eqref{eq:411a}, i.e.,
\beq\label{eq:41a}
a_1(\tu_h,v_h)=\Ftilde(v_h) \quad\text{ for all } v_h\in \Vhp,
\eeq
and let $\tbu$ be the vector of nodal values of $\tu_h$. The definition of $\widetilde{f}$ then implies that \eqref{eq:41a} is equivalent to $\Amato \tbu = \Smat_A\,\fvec$. 

Similar to the proof of the bounds on $\AmatoI\Mmatn$ in Lemma \ref{lem:keylemma1},
using the triangle inequality, the bound \eqref{eq:bound4} from Condition \ref{cond:2}, the bound \eqref{eq:bound2} from Lemma \ref{lem:H1}, the bound \eqref{eq:keyAboundfinal}, and the definition \eqref{eq:C1} of $C_1$,
we find
\begin{align}\nonumber 
\N{\tu_h}_{\HokD} &\leq
\N{\tu-\tu_h}_{\HokD} + \N{\tu}_{\HokD},\nonumber \\ \nonumber
& \leq \left[ C^{(1)}_{\rm FEM2} k + 
\frac{1}{\min\{{A_1}_{\min},{n_1}_{\min}\}}\left( 1 + 2 C^{(1)}_{\rm bound}{n_1}_{\max} k  \right) 
\right]\big\|\Ftilde\big\|_{(\HokD)^*},\\
&\leq C_1 \, k\, 
\Cinvs \NLqDRRdtd{A} h^{d\mleft(\frac1{s} - \half\mright)} \big\|\grad\ftilde\big\|_{L^2(D)},
\label{eq:mainevent2}\\
&\leq C_1 \, k\, 
\Cinvs\NLqDRRdtd{A} h^{d\mleft(\frac1{s} - \half\mright)}
 \big\|\widetilde{f}\big\|_{\HokD},\nonumber
\end{align}
and the bound on $\|(\Amato)^{-1}\Smat_A\|_{\Dmat_k}$ in \eqref{eq:keybound22} follows since $1/q = 1/2-1/s$.

To prove the bound on  $\|(\Amato)^{-1}\Smat_A\|_{2}$ in \eqref{eq:keybound2a2}, we use the bounds 
\beqs
m_- h^{d/2} k \N{\tbu}_2 \leq k \N{\widetilde{u}_h}_{\LtD} \leq \N{\widetilde{u}_h}_{\HokD}
\,\tand\,
\big\|\nabla \widetilde{f}\big\|_{\LtD} \leq s_+ h^{d/2-1}\N{\fvec}_2,
\eeqs
on either side of the inequality \eqref{eq:mainevent2}, with these bounds coming from \eqref{eq:normequiv1} and \eqref{eq:normequiv2} respectively. 
\epf

\bre[Link to the results of \cite{GaGrSp:15}]\label{rem:GaGrSp}
A result analogous to the Euclidean-norm bounds in Lemma \ref{lem:mainingredient} was proved in \cite[Theorem 1.4]{GaGrSp:15} for the case that $A_1= A_2= I$, $n_2= 1$, and $n_1 = 1 + \ri\eps/k^2$, with the `absorption parameter' or `shift' $\eps$ satisfying $0<\eps\lesssim k^2$. This is equivalent to approximating the discrete inverse of $\Delta + k^2$ by the discrete inverse of $\Delta +k^2 + \ri \eps$.
 The motivation for proving this result was that the so-called `shifted Laplacian preconditioning' of the Helmholtz equation (introduced in \cite{ErVuOo:04}) is based on preconditioning, with these choices of parameters, $\Amatt$ with an approximation of $\Amato$. Similar to the proofs of Theorems  \ref{thm:1new} and \ref{thm:1anew} in \S\ref{sec:mainproofs}, bounds on $\|\Imat -  (\Amato)^{-1}\Amatt \|_2$ and 
$\|\Imat - \Amatt  (\Amato)^{-1}\|_2$
 then give upper bounds on how large the `shift' $\eps$ can be for GMRES for $\AmatoI\Amatt$ to converge in a $k$-independent number of iterations in the case when the action of $(\Amato)^{-1}$ is computed exactly.

The main differences between \cite{GaGrSp:15} and this work are that (i)  \cite{GaGrSp:15} focuses only on the case when $T=\ri k$,
(ii) the theory in \cite{GaGrSp:15} focuses on the particular case that $\Dm$ is star-shaped with respect to a ball, finding a $k$- and $\eps$-explicit expression for $C^{(1)}_{\rm bound}$ in this case using Morawetz identities, (iii) the proof of \cite[Theorem 1.4]{GaGrSp:15} requires bounds on 
$(\Amato)^{-1}\Mmat_{n}$, $\Mmat_{n}(\Amato)^{-1}$ (analogous to the bounds in Lemma \ref{lem:keylemma1}), $(\Amato)^{-1}\Nmat$, and $\Nmat(\Amato)^{-1}$, but \emph{not} on 
$(\Amato)^{-1}\Smat_{A}$ and $\Smat_{A}(\Amato)^{-1}$, (iv) \cite[Theorem 1.4]{GaGrSp:15} holds under the assumption that the finite-element approximation is quasi-optimal, whereas here we prove results under Condition \ref{cond:2}, which holds  under less restrictive conditions on $h$ and $p$ than those for quasioptimality (see \S\ref{sec:cond2hold}), and (v) \cite{GaGrSp:15} only proves bounds
on $\Imat -  (\Amato)^{-1}\Amatt$ and  $\Imat - \Amatt  (\Amato)^{-1}$
 in the $\|\cdot\|_2$ norm whereas here we also prove bounds in weighted norms.
\ere

\subsection{Proofs of Theorems \ref{thm:Linfty1} and \ref{thm:Linfty2} from Theorems \ref{thm:1new} and \ref{thm:1anew}}
\label{sec:proofLinfty1}

Theorems \ref{thm:Linfty1} and \ref{thm:Linfty2} follow from 
repeating the proofs of 
Theorems \ref{thm:1new} and \ref{thm:1anew} but now, instead of using H\"older's inequality \eqref{eq:keynbound2}, using the simpler inequality 
\begin{equation*}
\big\|n\widetilde{f}\big\|_{L^2(D)} \leq \big\| n\big\|_{L^\infty(D)}\big||\widetilde{f}\big\|_{L^2(D)}.
\end{equation*}
The result is that the inverse inequality \eqref{eq:inverses} is no longer needed in \eqref{eq:mainevent1Euan} to 
obtain $\|\ftilde\|_{L^2(D)}$ from $\|\ftilde\|_{L^s(D)}$. 
Then Lemmas \ref{lem:keylemma1}, \ref{lem:keylemma2}, and \ref{lem:mainingredient}, 
and Theorems \ref{thm:1new} and \ref{thm:1anew} all 
hold without the factors $h^{-d/q}$ on the right-hand sides of the inequalities,
and with the constants $\widetilde{C}_j, j=1,2,$ replaced by $C_j, j=1,2$.
This is equivalent to these results 
holding with $q$ set formally as $\infty$, since 
$\Cinvs=1$ when $s=2$ (which occurs when $q=\infty$), as noted immediately after \eqref{eq:inverses}.

The end result is that the condition \eqref{eq:condalt} in Theorem \ref{thm:1new} is replaced by the condition \eqref{eq:cond} in Theorem \ref{thm:Linfty1} and similarly \eqref{eq:condaalt} in Theorem \ref{thm:1anew} is replaced by \eqref{eq:conda} in Theorem \ref{thm:Linfty2}.

\subsection{Proof of Lemma \ref{lem:standardvsweighted}}

Let $\br^\ell$ denote the $\ell$th residual of weighted GMRES in the norm $\|\cdot\|_{\Dmat_k}$
(as in the proof of Theorem \ref{thm:1new} in \S\ref{sec:mainproofs}), and let $\br^\ell_s$ denote the $\ell$th residual of standard GMRES.
By the minimal residual property \eqref{eq:minres},
\beqs
\N{\br^\ell_s}_2 \leq \N{\br^\ell}_2 \quad\tfa \ell\in \mathbb{N}.
\eeqs
By assumption $\br^0_s=\br^0$, and so
\beq\label{eq:Shihua1}
\frac{
\N{\br^\ell_s}_2}{
\N{\br^0_s}_2
} \leq 
\frac{\N{\br^\ell}_2}{
\N{\br^0}_2
} \quad\tfa \ell \in \mathbb{N}.
\eeq
By Lemma \ref{lem:normequiv}, the norm property \eqref{eq:Dk2}, and the fact that $hk \leq Ck^{-\delta}$,
there exists $\widetilde{C}>0$, dependent only on $s_+, m_+, k_0, C,$ and $\delta$, such that 
\beq\label{eq:Shihua2}
m_1 k h^{d/2} \N{\bv}_2 \leq \N{\bv}_{\Dmat_k} \leq \widetilde{C} h^{d/2-1}  \N{\bv}_2 \quad \tfa k\geq k_0.
\eeq
Therefore, by the combination of \eqref{eq:Shihua1}, \eqref{eq:Shihua2}, and the fact that $hk=Ck^{-\delta}$,
\beq\label{eq:Shihua3}
\frac{
\N{\br^\ell_s}_2}{
\N{\br^0_s}_2
} \leq 
\frac{\widetilde{C}}{m_- (hk)} 
\frac{\N{\br^\ell}_{\Dmat_k}}{
\N{\br^0}_{\Dmat_k}
}
\leq \frac{\widetilde{C}}{m_- C} k^\delta 
\frac{\N{\br^\ell}_{\Dmat_k}}{
\N{\br^0}_{\Dmat_k} 
}
 \quad \tfa \ell \in \mathbb{N} \,\tand\, k\geq k_0.
\eeq
Since \eqref{eq:cond} holds, the proof of Theorem \ref{thm:Linfty1} (in particular \eqref{eq:GMREScm})
shows that 
\beq\label{eq:Shihua4}
\big\| \rvec^\ell\big\|_{\Dmat_k} \leq c^\ell \big\|\rvec^0\big\|_{\Dmat_k} \quad\tfa \ell \in \mathbb{N}, 
\eeq
with $c:= 2\sqrt{2}/3<1$. Therefore, combining \eqref{eq:Shihua3} and \eqref{eq:Shihua4}, we find that
\beqs
\frac{
\N{\br^{\ell+n}_s}_2}{
\N{\br^0_s}_2
} 
\leq \frac{\widetilde{C}}{m_- C} k^\delta 
\frac{\N{\br^{\ell+n}}_{\Dmat_k}}{
\N{\br^0}_{\Dmat_k} 
}
\leq \frac{\widetilde{C}}{m_- C} k^\delta 
c^{\ell+n}
 \quad \tfa \ell, n \in \mathbb{N} \,\tand\, k\geq k_0.
\eeqs
Therefore, if $n$ is such that 
\beq\label{eq:conditiononn}
\frac{\widetilde{C}}{m_- C} k^\delta c^n \leq 1,
\eeq
and $\ell \geq \log (1/\eps)/\log(1/c)$, 
then 
\beqs
\frac{
\N{\br^{\ell+n}_s}_2}{
\N{\br^0_s}_2
} \leq \eps.
\eeqs
Arguing exactly as in \eqref{eq:interim1} in the proof of Theorem \ref{thm:Linfty1}, we see that the result of the lemma, i.e., the bound 
\eqref{eq:compare1b}, holds if $m \geq \log (1/\eps)/\log(1/c) +n$ where $n$ is such that \eqref{eq:conditiononn} holds.
The lower bound \eqref{eq:compare1a} on $m$ then follows with $\mathcal{C}:= \log(\widetilde{C}/ (m_- C))$, since by \eqref{eq:conditiononn} $n \geq (\mathcal{C} + \delta\log k)/\log(1/c)$. 

\section{Proofs of Theorem \ref{thm:2} and Lemma \ref{lem:sharp}}\label{sec:proofsPDE}

\bpf[Proof of Theorem \ref{thm:2}]
Since $a_1(u_1, v) = F(v)$ and 
$a_2(u_2, v) = F(v)$, we have
\beq\label{eq:vp1}
a_1(u_1-u_2,v) = \widetilde{F}(v) \quad\tfa v\in H^1_{0,D}({D})
\eeq
where
\beqs
 \widetilde{F}(v):= \int_{{D}} \left((A_2-A_1) \nabla u_2\right) \cdot\overline{\nabla v} + k^2 ({n_1}-{n_2}) u_2\overline{v}.
\eeqs
Now, by the Cauchy-Schwarz inequality and the definition \eqref{eq:1knorm} of the norm $\|\cdot\|_{\HoDkk}$, we have that
\begin{align*}
| \widetilde{F}(v)| &\leq \big\|A_1-A_2\big\|_{L^\infty(D; {\rm op})} \big\|\nabla u_2\big\|_{L^2({D})}
\N{\nabla v}_{L^2({D})} 
\\& \hspace{5cm}+ k^2 
\big\|{n_1}-{n_2}\big\|_{L^\infty(D;\Rea)} \big\| u_2\big\|_{L^2({D})}
\N{v}_{L^2({D})}\\
&\leq\max\Big\{\big\|A_1-A_2\big\|_{L^\infty(D; {\rm op})}\,,\, \big\|{n_1}-{n_2}\big\|_{L^\infty(D;\Rea)}\Big\}
\big\| u_2\big\|_{\HoDkk} \N{v}_{\HoDkk}.
\end{align*}
and thus, by the definition \eqref{eq:dualnorm} of the norm $\|\cdot\|_{(\HoDkk)*}$,
\beqs
\big\|\widetilde{F}\big\|_{(\HoDkk)*}\leq \max\Big\{\big\|A_1-A_2\big\|_{L^\infty(D; {\rm op})}\,,\, \big\|{n_1}-{n_2}\big\|_{L^\infty(D;\Rea)}\Big\}
\big\| u_2\big\|_{\HoDkk}.
\eeqs
Since Condition \ref{cond:1} holds, we can then apply the result of Lemma \ref{lem:H1}, i.e.~the bound \eqref{eq:bound2}, to the solution of the variational problem \eqref{eq:vp1}  to find that 
\begin{align*}
\frac{\big\| u_1 - u_2\big\|_{\HoDkk}}
{\big\| u_2\big\|_{\HoDkk}
}
 \leq 
\,&\frac{1}{\min\big\{{A_1}_{\min},{{n_1}}_{\min}\big\}}\left( 1 + 2 C^{(1)}_{\rm bound}{{n_1}}_{\max}  k\right)
\\
&\quad\times \left(\max\Big\{\big\|A_1-A_2\big\|_{L^\infty(D; {\rm op})}\,,\, \big\|{n_1}-{n_2}\big\|_{L^\infty(D;\Rea)}\Big\}\right),
\end{align*}
and then the result \eqref{eq:PDEbound} follows with 
\beq\label{eq:C3}
C_3:= \frac{1}{\min\big\{{A_1}_{\min},{{n_1}}_{\min}\big\}}\left( \frac{1}{k_0} + 2 C^{(1)}_{\rm bound}{{n_1}}_{\max}  \right).
\eeq
\epf

\bpf[Proof of Lemma \ref{lem:sharp}]
We actually prove the stronger result that, given any function $c(k)$ such that $c(k)>0$ for all $k\geq k_0$, there exist 
$f, {n_1}, {n_2}$, with ${n_1}\not= {n_2}$ and 
\beqs
\big\|{n_1}-{n_2}\big\|_{L^\infty(D;\Rea)} \sim c(k),
\eeqs
such that the corresponding solutions $u_1$ and $u_2$ of the exterior Dirichlet problem with $A_1 = A_2= I$ exist, are unique, and satisfy \eqref{eq:sharp1}. 

The heart of the proof of \eqref{eq:sharp1} is the equation
\beq\label{eq:obs1}
(\Delta + k^2) \big(\re^{\ri k r}\chi(r)\big) =  \re^{\ri k r}\left(\ri k \frac{d-1}{r} \chi(r) + 2 \ri k \diff{\chi}{r}(r) + \Delta \chi(r)\right)=: -\widetilde{f}(r),
\eeq
where $\chi(r)$ is chosen to have $\supp \chi \subset D$.
This example proves the sharpness of the nontrapping resolvent estimate \eqref{eq:bound1} as $k\tendi$, since both the $L^2(D)$ norm of $\widetilde{f}$ and the $\HoDkk$ norm of $\re^{\ri kr}\chi(r)$ are proportional to $k$, and hence each to other (see, e.g., \cite[Lemma 3.10]{ChMo:08},  \cite[Lemma 4.12]{Sp:14}).

The overall idea of the proof of \eqref{eq:sharp1} is to set things up so that $(u_1-u_2)(\bx) = \re^{\ri k r}\chi(r)$, the rationale being that \eqref{eq:obs1} proves the sharpness of \eqref{eq:bound1}, and \eqref{eq:bound1} and its corollary \eqref{eq:bound2} (applied to $u_1-u_2$) are the main ingredients in the proof of Theorem \ref{thm:2}.

With the choices $A_j:=I$, $j=1,2,$ and ${n_1}:=1$, the variational problem \eqref{eq:vp1} implies that 
\beq\label{eq:obs2}
(\Delta +k^2)(u_1 - u_2) = k^2 \big(n_2-1\big)u_2.
\eeq
We now fix ${n_2}:= 1 + c(k)\widetilde{\chi}$ with  $\widetilde{\chi}\in C^{\infty}({D})$ chosen so that $\widetilde{\chi}= \widetilde{\chi}(r)$,
 $\supp \, \widetilde{\chi} \subset{D}$, $\widetilde{\chi}$ is positive everywhere in the interior of its support, and $\|\widetilde{\chi}\|_{L^\infty(D;\mathbb{R})}=1$; observe that this last condition implies that  $\|{n_1}- {n_2}\|_{L^\infty(D;\Rea)} = c(k)$.
 As above, let $\chi=\chi(r)$ with $\chi \in C^{\infty}({D})$ and
$\supp \,\chi\subset{D}$. Then, with $\widetilde{f}$ defined by \eqref{eq:obs1}, we formally define $u_2$ and $f$ by 
\beq\label{eq:obs3}
u_2(\bx):= -\frac{1}{k^2 c(k)}\frac{\widetilde{f}(r)}{\widetilde{\chi}(r)}, \quad\tand \quad f(\bx):= -\big(\Delta +k^2 {n_2}(\bx)\big) u_2(\bx),
\eeq
and we define $u_1$ to be the solution of the TEDP (in the sense of Definition \ref{def:TEDP}) with $n_1=1$, $g=0$, and $f$ defined as in \eqref{eq:obs3}. Since $\widetilde{\chi}(r)$ has compact support in $D$, we need to tie both the support of $\widetilde{\chi}$ and how fast $\widetilde{\chi}$ vanishes near the boundary of its support to the definition of $\chi$ for both $u_2$ and $f$ in \eqref{eq:obs3} to be well defined.
We ignore this issue for the moment and show that these definitions achieve our goal of having
\beq\label{eq:obs4}
u_1(\bx)- u_2(\bx) = \re^{\ri k x_1}\chi(r).
\eeq
First observe that
since $\supp \,\widetilde{f}$ is a compact subset of ${D}$, so is 
$\supp \,u_2$. Therefore $u_2$ is the solution of the TEDP 
with $g=0$, $f$ defined in \eqref{eq:obs3}, and coefficient ${n_2}:=1 + c(k)\widetilde{\chi}$.
Combining  \eqref{eq:obs2} with the definitions of $n_2$ and $u_2$, we see that 
\beq\label{eq:obs5}
(\Delta +k^2)(u_1-u_2) = -\widetilde{f}.
\eeq
Since the solution of the TEDP is unique (see Remark \ref{rem:existence}), \eqref{eq:obs1} and \eqref{eq:obs5} imply that \eqref{eq:obs4} holds. We therefore have that
\beqs
\big\|u_1-u_2\big\|_{L^2({D})} \sim 1
\quad \tand \quad
\big\|u_1-u_2\big\|_{\HoDkk} \sim k.
\eeqs
Furthermore, the definitions \eqref{eq:obs3} and \eqref{eq:obs1} of $u_2$  and $\widetilde{f}$, respectively, imply that
\beqs
\big\| u_2\big\|_{L^2({D})} \sim \frac{1}{k\, c(k)} \quad\tand \quad 
\big\| u_2\big\|_{\HoDkk} \sim \frac{1}{c(k)},
\eeqs 
and, since $\|{n_1}- {n_2}\|_{L^\infty(D;\Rea)} = c(k)$, \eqref{eq:sharp1} holds.

Therefore, to complete the proof, we only need to show that there exists a choice of $\chi$ and $\widetilde{\chi}$ for which $u_2$ and $f$ defined by \eqref{eq:obs3} are 
in $H^{1}(\OR)$ and $L^2(\OR)$ respectively (in fact, we prove that they are in $W^{1,\infty}(\OR)$ and $L^\infty(\OR)$ respectively). 
Since $\chi$ and $\widetilde{\chi} \in C^\infty(\OR)$, 
and $\widetilde{\chi}$ is positive everywhere in the interior of its support,
the only issue is what happens at the boundary of $\supp\widetilde{\chi}$, where $u_2$ could be singular.

By assumption there exist
 $0<R_1<R_2$ such that $\overline{\Dm} \subset B_{R_1}\subset B_{R_2} \subset  \widetilde{D}$, i.e.~there exists an annulus between $D_-$ and $\partial \widetilde{D}$.
 Let $\supp \chi = B_{R_2}\setminus B_{R_1}$ and let $\chi$ vanish to order $m\geq 3$ at $r= R_1$ and $r=R_2$; i.e.~$\chi(r) \sim (r-R_1)^m$ as $r \rightarrow (R_1)^+$ and 
$\chi(r) \sim (R_2-r)^m$ as $r \rightarrow (R_2)^-$. The definition  \eqref{eq:obs1} of $\widetilde{f}$ then implies that $\widetilde{f}$ vanishes to order $m-2$. Let $\widetilde{\chi}(r)$ vanish to order $\ell\geq 1$ at $r= R_1$ and $r=R_2$. 
We now claim that if $m >\ell+4$, then $u_2\in W^{1,\infty}(\OR)$ and $f$ $\in L^\infty({D})$. Indeed,  
straightforward calculation from \eqref{eq:obs3} shows that  $u_2(r) \sim (r-R_1)^{m-\ell-2}$, $\nabla u_2(r) \sim (r-R_1)^{m-\ell-3}$, and $\Delta u_2(r) \sim (r-R_1)^{m-\ell-4}$ as $r \rightarrow (R_1)^+$, with analogous behaviour at $r=R_2$.
The assumption 
$m >\ell+4$ therefore implies that $u_2$, $\nabla u ^{(2)}$, and $\Delta u_2$ vanish (and hence are finite) at $r=R_1$ and $r=R_2$.
\epf

\bre[Why doesn't Lemma \ref{lem:sharp} cover the case $A_1\neq  A_2$?]
When $n_j:=1$, $j=1,2,$ $A_1:=I$, and $A_2:= I + c(k)\widetilde{\chi}$, the variational problem \eqref{eq:vp1} implies that 
\beqs
\Delta \big( u_1 - u_2\big) + k^2 \big( u_1 - u_2\big) = c(k)\nabla\cdot \big(\widetilde{\chi}\nabla u_2\big).
\eeqs
It is now much harder than in \eqref{eq:obs2} to set things up so that $ u_1(\bx) - u_2(\bx)=\re^{\ri kr}\chi(r)$ (so that one can then use \eqref{eq:obs1}).
\ere

\bre[The relationship between Q1 and \Qtwo]\label{rem:analogue}
In \S\ref{sec:sharpness}, we stated that \Qtwo is an analogue of Q1 on the PDE level; we now make this statement more precise.

Since $C_3$ defined by \eqref{eq:C3} is independent of $F$, the relative-error bound \eqref{eq:PDEbound} of Theorem \ref{thm:2} can be written as
\beq\label{eq:PDEbound2}
\sup_{F\in( H^1_{0,D}(D))^*}
 \frac{\big\|u_1-u_2\big\|_{\HokD}
}{
\N{u_2}_{\HokD}
}\leq C_3 \,k\, \max\set{\NLiDRRdtd{A_1-A_2}\,,\, \NLiDRR{n_1-n_2}}.
\eeq

Let $\cA: H^1_{0,D}(D)\rightarrow ( H^1_{0,D}(D))^*$ be the operator associated with the sesquilinear form $a(\cdot,\cdot)$ \eqref{eq:agen} so that the variational problem \eqref{eq:vp} is equivalent to the operator equation $\cA u = F$ in $( H^1_{0,D}(D))^*$ (for details of this equivalence, see, e.g., \cite[Lemma 2.1.38]{SaSc:11}). 
Then 
$\cA_j u_j= F$ 
and the assumptions in Theorem \ref{thm:2} 
imply that $\cA_1$ and $\cA_2$ are both invertible.
Using this fact, we see that the left-hand side of \eqref{eq:PDEbound2} equals
\begin{align*}
\sup_{F\in( H^1_{0,D}(D))^*}
 \frac{\big\|(\cA_1)^{-1}F- (\cA_2)^{-1}F\big\|_{\HokD}
}{
\N{(\cA_2)^{-1}F}_{\HokD}
}
&=
\sup_{F\in( H^1_{0,D}(D))^*}
 \frac{\big\|\big( I - (\cA_1)^{-1}\cA_2\big)(\cA_2)^{-1}F\big\|_{\HokD}
}{
\N{(\cA_2)^{-1}F}_{\HokD}
}\\
&= \Big\|  I - (\cA_1)^{-1}\cA_2\Big\|_{\HokD\rightarrow\HokD}
\end{align*}
so that \eqref{eq:PDEbound2} becomes 
\beq\label{eq:PDEbound3}
\Big\|  I - (\cA_1)^{-1}\cA_2\Big\|_{\HokD\rightarrow\HokD}\leq C_3 \,k\, \max\set{\NLiDRRdtd{A_1-A_2}\,,\, \NLiDRR{n_1-n_2}}.
\eeq
On the other hand, the main ingredient to proving Theorem \ref{thm:Linfty1} is the bound 
\beqs
\Big\|\Imat - \AmatoI\Amatt\Big\|_{\Dmat_k}\leq \Cotilde \,k \NLiDRRdtd{A_1-A_2} + \Cttilde  \,k  \NLiDRR{n_1-n_2},
\eeqs
which implies
\beq
\Big\|\Imat - \AmatoI\Amatt\Big\|_{\Dmat_k}\leq 2 \,k\, \max \left\{ \Cotilde  \NLiDRRdtd{A_1-A_2}\,,\, \Cttilde   \NLiDRR{n_1-n_2}\right\}. 
\label{eq:discretebound1}
\eeq
The left-hand side of \eqref{eq:discretebound1} is the discrete analogue of the left-hand side of \eqref{eq:PDEbound3}.
The right-hand sides of these bounds are (up to constants independent of $k$) identical, hence why we say that the condition \eqref{eq:sufficientlysmall} (i.e.~$k\|A_1-A_2\|_{L^\infty(D;{\rm op})}$ and $k\|n_1-n_2\|_{L^\infty(D;\Rea)}$ both sufficiently small) is the answer to both Q1 and \Qtwo.
\ere

\section*{Acknowledgements.}

We thank  Th\'eophile Chaumont-Frelet (INRIA, Nice), Stefan Sauter (Universit\"at Z\"urich), and Nilima Nigam (Simon Fraser University) for useful comments and discussions about this work at the conference MAFELAP 2019. We thank the referees for their constructive comments and insightful suggestions. Finally, we thank Ralf Hiptmair (ETH Z\"urich) and Robert Scheichl (Universit\"at Heidelberg) for useful comments on this work in the course of examining ORP's PhD thesis \cite{Pe:20}.

IGG acknowledges support from EPSRC grant  EP/S003975/1.
ORP acknowledges support by a scholarship from the EPSRC Centre for Doctoral Training in Statistical Applied Mathematics at Bath (SAMBa), under the EPSRC grant EP/L015684/1.
EAS acknowledges support from EPSRC grant EP/R005591/1.
This research made use of the Balena High Performance Computing (HPC) Service at the University of Bath.
    
\footnotesize{
\bibliographystyle{plain}

}
\end{document}